\documentclass[11pt]{article}

\usepackage{latexsym}
\usepackage{amssymb}
\usepackage{amsthm}
\usepackage{amscd}
\usepackage{amsmath}
\usepackage{mathrsfs}
\usepackage[all]{xy}
\usepackage{graphicx}
\usepackage{hyperref} 

\input xy \xyoption{frame}

\usepackage[usenames,dvipsnames]{color}

\theoremstyle{definition}

\newtheorem* {theorem*}{Theorem}

\newtheorem{theorem}{Theorem}[section]
\newtheorem{question}{Question}[section]
\theoremstyle{definition}
\newtheorem{observation}{Observation}[section]

\newtheorem* {example*}{Example}

\newtheorem{lemma}{Lemma}[section]
\theoremstyle{definition}

\theoremstyle{definition}

\newtheorem{proposition}{Proposition}[section]
\newtheorem{corollary}{Corollary}[section]
\newtheorem* {remark}{Remark}
\theoremstyle{definition}

\theoremstyle{definition}

\theoremstyle{definition}

\theoremstyle{definition}

\xyoption{dvips}

\usepackage{fullpage}

\numberwithin{equation}{section}

\newcommand{\One}{{1\hspace{-.14cm} 1}}

\def\sG{\mathcal{G}}
\def\cM{\mathcal{M}}
\def\Fix{\mathrm{Fix}}

\def\modu{\ (\mathrm{mod}\ }

\def\({\left(}
\def\){\right)}

\newcommand{\sgn}{\mathrm{sgn}}

                \newcommand{\cX}{\mathcal{X}}
    \newcommand{\sP}{\Lambda}      \newcommand{\CC}{\mathbb{C}}  \newcommand{\RR}{\mathbb{R}} \newcommand{\QQ}{\mathbb{Q}}    
    \newcommand{\cC}{\mathcal{C}}
\newcommand{\cR}{\mathcal{R}}   \newcommand{\cI}{\mathcal{I}}

    \def\ZZ{\mathbb{Z}} \def\Aut{\mathrm{Aut}}   \def\Ind{\mathrm{Ind}} \def\GL{\mathrm{GL}}   \def\Res{\mathrm{Res}}    \def\spanning{\textnormal{-span}}   
\def\Irr{\mathrm{Irr}}  \def\wt{\mathrm{wt}}

\def\sign{\mathrm{sign}}

\def\Inv{\mathrm{Inv}}
\def\Pair{\mathrm{Pair}}
\def\Inn{\mathrm{Inn}}
\def\Ad{\mathrm{Ad}}

\def\wt{\widetilde}
\def\sP{\mathscr{P}}

\def\cV{\mathcal{V}}

\def\Out{\mathrm{Out}}
\def\cyc{\ }

\newcommand{\leftexp}[2]{{\vphantom{#2}}^{#1}{#2}}
\newcommand{\ba}{\begin{aligned}}
\newcommand{\ea}{\end{aligned}}
\newcommand{\barr}{\begin{array}}
\newcommand{\earr}{\end{array}}
\newcommand{\be}{\begin{equation}}
\newcommand{\ee}{\end{equation}}
\makeatletter
\renewcommand{\@makefnmark}{\mbox{\textsuperscript{}}}
\makeatother

\allowdisplaybreaks[1]

\UseCrayolaColors

\begin{document}
\title{Automorphisms and generalized involution models of finite complex reflection groups}
\author{Eric Marberg\footnote{This research was conducted with government support under
the Department of Defense, Air Force Office of Scientific Research, National Defense Science
and Engineering Graduate (NDSEG) Fellowship, 32 CFR 168a.} \\ Department of Mathematics \\ Massachusetts Institute of Technology, United States \\ \tt{emarberg@math.mit.edu}}
\date{}

\maketitle

\begin{abstract}
    We prove that a finite complex reflection group has a generalized involution model, as defined by Bump and Ginzburg,  if and only if each of its irreducible factors is either $G(r,p,n)$ with $\gcd(p,n)=1$; $G(r,p,2)$ with $r/p$ odd; or $G_{23}$, the Coxeter group of type $H_3$.  We additionally provide explicit formulas for all automorphisms of $G(r,p,n)$, and construct new Gelfand models for the groups $G(r,p,n)$ with $\gcd(p,n)=1$.   
\end{abstract}

\section{Introduction}

A \emph{model} for a finite group $G$ is a set $\{ \lambda_i : H_i \to \CC\}$ of linear characters  of subgroups of $G$ such that
$ \sum_i \Ind_{H_i}^G(\lambda_i)$ is the multiplicity-free sum of all irreducible characters of $G$.  If all of the subgroups $H_i$ are centralizers of involutions, with each conjugacy class of involutions contributing exactly one subgroup, then we say that a model is an \emph{involution model}.  Bump and Ginzburg introduced in \cite{BG2004} the notion of a \emph{generalized involution model} as a natural extension of this concept.

Generalized involution models are parametrized by automorphisms which are involutions.    Write the action of $\tau \in \Aut(g)$ on $g \in G$ by $\leftexp{\tau}g$.  
For each automorphism $\tau \in \Aut(G)$ with $\tau^2=1$, 
the group $G$ acts on the set of \emph{generalized involutions} 
 \[\cI_{G,\tau} \overset{\mathrm{def}} = \{ \omega \in G : \omega\cdot \leftexp{\tau}{\omega} = 1\}\] by the $\tau$-twisted conjugation $g : \omega \mapsto  g\cdot \omega\cdot \leftexp{\tau}{g}^{-1}.$  Let 
 \[ C_{G,\tau}(\omega) = \{ g \in G : g\cdot \omega \cdot \leftexp{\tau}g^{-1}=\omega\}\] denote the stabilizer of $\omega \in \cI_{G,\tau}$ in $G$ under this action.  We refer to the orbit of $\omega \in \cI_{G,\tau}$ as a \emph{twisted conjugacy class} and to $C_{G,\tau}(\omega)$ as its \emph{twisted centralizer}.  

A \emph{generalized involution model} for $G$ with respect to $\tau$ is a model $\cM$ with an injective map $\iota: \cM\to \cI_{G,\tau}$ such that the following hold:
\begin{enumerate}
\item[(a)] Each $\lambda \in \cM$ is a linear character of the $\tau$-twisted centralizer in $G$ of $\iota(\lambda) \in \cI_{G,\tau}$.
\item[(b)] The image of $\iota$ contains exactly one element from each $\tau$-twisted conjugacy class in $\cI_{G,\tau}$.
\end{enumerate}
This definition differs slightly from the one originally given in \cite{BG2004}, but one can show that it is equivalent, in the sense that the same models are classified as generalized involution models.

Call a representation of $G$ equivalent to the multiplicity-free sum of all of the group's irreducible representations a \emph{Gelfand model}.   One reason why generalized involution models are interesting mathematical objects is that they often ``explain'' the construction of elegant Gelfand models.  As a motivating example, Adin, Postnikov, and Roichman recently described a combinatorial Gelfand model for the wreath product $\ZZ_r \wr S_n$ \cite{APR2007,APR2008}.  In  \cite{M} we discuss how their construction arises naturally from the merging together of generalized involution models for  abelian and symmetric groups.
Moreover, the existence of a generalized involution model provides a combinatorial interpretation of several properties of a group's character table, such as its row sums.  

In comparison to their classical predecessor, generalized involution models are more flexible and, contrary to appearances, in practice often not much more difficult to classify.  Past researchers have investigated classification questions concerning involution models.  Inglis, Richardson, and Saxl  identified an involution model for the symmetric group in \cite{IRS91}.  Not long after, Baddeley classified in his doctoral thesis which irreducible Weyl groups have involution models: in particular, only those of types $D_{2n}$ ($n>1$), $F_4$, $E_6$, $E_7$, and $E_8$ do not have involution models \cite{B91}.  Vinroot recently extended this classification to finite Coxeter groups, in particular showing that the Coxeter groups of type $I_2(n)$ and $H_3$ each have an involution model while the group of type $H_4$ does not \cite{V}.

Investigating which complex reflection groups have involution models appears at first glance a logical extension of this body of work.  Involution models, however, exist only for groups whose irreducible representations are all realizable over $\RR$.
Thus, the right question to ask$-$and the question answered in this work$-$is which complex reflection groups have \emph{generalized} involution models.  Specifically, we will prove the following theorem.

\begin{theorem*}
A finite complex reflection group has a generalized involution model if and only if each of its irreducible factors is one of the following:
\begin{enumerate}
\item[(i)] $G(r,p,n)$ with $\gcd(p,n)=1$.
\item[(ii)] $G(r,p,2)$ with $r/p$ odd.
\item[(iii)] $G_{23}$, the Coxeter group of type $H_3$.
\end{enumerate}
\end{theorem*}

In the course of proving this result, we construct new Gelfand models for the irreducible complex reflection groups $G(r,p,n)$ with $\gcd(p,n)=1$, different from the ones described  by Caselli in the recent paper \cite{C2009}.  We additionally provide explicit formulas for the automorphisms of the complex reflection groups $G(r,p,n)$, deriving among other things an expression for the order of $\Aut\(G(r,p,n)\)$.  We intend this analysis to complement the general characterizations of the automorphisms of complex reflection groups given in the recent work \cite{MM} and the classification of the subgroup of reflection-preserving automorphisms in \cite{ShiWang}.


\section{Preliminaries}\label{prelim}

In this section we review some useful facts concerning complex reflection groups and their irreducible characters, and generalized involutions models.

\subsection{Complex Reflection Groups}

A \emph{pseudo-reflection} is an automorphism of a complex vector space which fixes every point in some hyperplane.  A \emph{complex reflection group} $G$ is a group generated by a set of pseudo-reflections of some finite-dimensional complex vector space $V$.  If no proper, nonzero subspaces of $V$ are $G$-invariant, then we say that $G$ is irreducible.  

In most situations, to study complex reflection groups it suffices to study the irreducible groups, since every finite complex reflection group decomposes as a direct product of irreducible complex reflection groups.  The finite irreducible complex reflection groups were identified through the work of a number of mathematicians in the nineteenth and first half of the twentieth century.  Shephard and Todd completed this classification in their seminal paper \cite{ST}; a useful modern treatment of this material appears in \cite{LT}.

The finite irreducible groups include one infinite family $G(r,p,n)$ and thirty-four exceptional groups labeled $G_4,\dots G_{37}$.  Presentations for the exceptions as abstract groups appear in \cite{BMR}.
We can describe the infinite series of groups $G(r,p,n)$ more concretely.  Let $r,p,n$ be positive integers with $p$ dividing $r$. 
As a subgroup of $\GL_n(\CC)$, the group $G(r,p,n)$ can be realized as the set of generalized permutation matrices whose nonzero entries are complex $r$th roots of unity, such that the product of the nonzero entries in any matrix is an $(r/p)$th root of unity.  This group acts irreducibly on $\CC^{n}$ when $r>1$ and $(r,p,n)\neq (2,2,2)$, and on the codimension 1 subspace of $\CC^n$ consisting of vectors whose coordinates sum to zero when $r=1$ and $n>1$.  

For our purposes, it will typically be more convenient to view $G(r,p,n)$ as a subgroup of the wreath product of a cyclic group with a symmetric group.  
The notation we will use for this goes as follows.  Given integers $a,b$, we write $[a,b]$ for the set $ \{ i \in \ZZ : a\leq i \leq b\}$. Fix positive integers $r,n$.  Let $\ZZ_r$ denote the cyclic group of order $r$, viewed as the additive group of integers $[0,r-1]$  modulo $r$, and let $S_n$ denote the symmetric group of permutations on $[1,n]$.  
The group $S_n$ acts on $(\ZZ_r)^n$ by permuting coordinates; we denote this action by
\[ \pi(x) \overset{\mathrm{def}} = \( x_{\pi^{-1}(1)}, \dots, x_{\pi^{-1}(n)}\),\qquad\text{for }\pi \in S_n,\ x=(x_1,\dots,x_n) \in (\ZZ_r)^n.\]  
The wreath product $\ZZ_r \wr S_n$ is the group given by the set of pairs $(x,\pi)$ with $x \in (\ZZ_r)^n$ and $\pi \in S_n$, with multiplication defined by
\[ (x,\pi)(y,\sigma) = (\sigma^{-1}(x)+  y, \pi \sigma),\qquad\text{for }x,y \in (\ZZ_r)^n,\ \pi,\sigma \in S_n.\]  Throughout, we identify  $(\ZZ_r)^n$ and $S_n$ with the subgroups $\{ (x,1) :x \in (\ZZ_r)^n\}$ and $\{(1,\pi) : \pi \in S_n\}$ in $\ZZ_r \wr S_n$.

Given an element $g = (x,\pi) \in \ZZ_r \wr S_n$, we define $|g| \in S_n$ and  $z_g : [1,n] \to \ZZ_r $ and  $\overline g, g^T \in \ZZ_r \wr S_n$ by 
\be\label{book} |g| = \pi,\qquad z_g(i) = x_i,\qquad \overline g = (-x,\pi),\qquad\text{and}\qquad g^T=\overline{g^{-1}} = \(\pi(x) , \pi^{-1} \).\ee  We call $g^T$ the \emph{transpose} of $g$, and say that $g$ is \emph{symmetric} if $g= g^T$.  After fixing a primitive $r$th root of unity $\zeta_r$, it makes sense to view each element $(x,\pi) \in \ZZ_r \wr S_n$ as the $n\times n$ matrix with $\zeta^{x_i}$ in the position $(\pi(i), i)$ for $i=1,\dots,n$ and zeros in all other positions.   If we identify $g$ with a generalized permutation matrix in this way, then $g^T$ corresponds to the usual matrix transpose of $g$, and $\overline g$ is the complex conjugate of $g$.   In particular, the map $g\mapsto g^T$ is an anti-automorphism and $g\mapsto \overline g=(g^{-1})^T$ is an automorphism, which we refer to as the \emph{inverse transpose automorphism}.

To define $G(r,p,n)$ we make use of the natural homomorphism $\Delta : (\ZZ_r)^n \to \ZZ_r$ given by 
\[ \Delta(x) = x_1 + x_2 + \dots + x_n,\qquad\text{for }x \in (\ZZ_r)^n.\]  This map extends to a homomorphism $\Delta:
 \ZZ_r \wr S_n \to \ZZ_r$  by the formula $\Delta(x,\pi) =\Delta(x)$ for $x \in (\ZZ_r)^n$ and $\pi \in S_n$.  Observe that $\Delta(\overline g)=\Delta(g^{-1}) = -\Delta(g)$ and $\Delta(g) = \Delta(g^T)$ for $g \in \ZZ_r \wr S_n$.

Given positive integers $r,p,n$,  we define the {complex reflection group} $G(r,p,n)$ as the normal subgroup of $\ZZ_r \wr S_n$ given by 
\[ G(r,p,n) 
 = \left\{ g \in \ZZ_r \wr S_n : \Delta(g) \in p\ZZ_r \right\}, \] where $p\ZZ_r$ denotes the subgroup $\{ 0,d ,2d, \dots, r-d\}\subset \ZZ_r$ generated by the greatest common divisor $d = \gcd(r,p)$.  In particular, if $d=r$ then $p\ZZ_r = \{0\}$.  To avoid redundancy in this definition, we henceforth require that $p$ divide $r$.   
 
 We observe that the wreath product $\ZZ_r \wr S_n$ is just $G(r,1,n)$.  
Likewise, every finite Coxeter group is a finite complex reflection group.  The Coxeter groups of type $A_n$, $B_n$, $D_n$, $G_2$, and $I_2(n)$ appear within the infinite series as  $G(1,1,n+1)$, $G(2,1,n)$, $G(2,2,n)$, $G(6,6,2)$, and $G(n,n,2)$  respectively.  The remaining finite Coxeter groups of type $H_3$, $F_4$, $H_4$, $E_6$, $E_7$, and $E_8$ appear as the exceptional groups $G_{23}$, $G_{28}$, $G_{30}$, $G_{35}$, $G_{36}$, and $G_{37}$, respectively.

\subsection{Irreducible Characters of $G(r,p,n)$}

To understand the irreducible characters of $G(r,p,n)$ we begin from a more general standpoint.  
Throughout, let $\Irr(G)$ denote the set of irreducible characters of a finite group $G$.  Now consider a finite group $G$ with a normal subgroup $H$ such that $G/H$ is cyclic.  Let $C \cong G/H$ denote the cyclic group of linear characters $\gamma$ of $G$ with $\ker \gamma \supset H$.  Then the tensor product $\otimes$ defines an action of $C$ on the irreducible characters of $G$, and we can say the following:
\begin{enumerate}
\item[(i)] Each irreducible character $\chi$ of $G$ restricts to a multiplicity-free sum of $k$ irreducible characters of $H$, where $k$ is the order of the stabilizer subgroup $\{\gamma \in C :  \gamma \otimes \chi = \chi \}$.  Furthermore, each irreducible character of $H$ is a constituent of some such restriction.

\item[(ii)] If $\chi,\psi$ are irreducible characters of $G$, then the following are equivalent:
\begin{enumerate}
\item[(a)] $\Res_H^G(\chi)$ and $\Res_H^G(\psi)$ share a common irreducible constituent.
\item[(b)] $\Res_H^G(\chi)=\Res_H^G(\psi)$.
\item[(c)] $\chi = \gamma \otimes \psi$ for some $\gamma \in C$.
\end{enumerate}
\end{enumerate}
These statements follow from Clifford theory; see \cite[Section 6]{S1989} for proofs.  

Specializing to the case at hand, we fix $r,p,n$ with $p$ dividing $r$ and let $G = G(r,1,n)$ and $H = G(r,p,n)$.  Observe that  $H$ is a normal subgroup of $G$ whose corresponding quotient group is cyclic and of order $p$.  
The irreducible characters of $G$ can be described as follows. 
For each $i \in [0,r-1]$, let $\psi_i : \ZZ_r \to \CC$ denote the irreducible character with the formula
\[ \psi_i(x) = \zeta_r^{ix},\qquad\text{for $x\in \ZZ_r$, where $\zeta_r$ is a fixed primitive $r$th root unity.}\]  Here we identify each $x \in \ZZ_r$ with an integer in $[0,r-1]$.  Then the trivial character $\One = \psi_0,$ and $\One, \psi_1,\psi_2,\dots,\psi_{r-1}$ are the distinct elements of $\Irr(\ZZ_r)$.  
Additionally, given a partition $\lambda$ of $n$, let $\chi^\lambda$ denote the corresponding irreducible character of $S_n$.  Let \[ \ba \sP &= \text{the set of all partitions of nonnegative integers,} \\
 \sP_r(n) &=\text{the set of $r$-tuples $\theta = (\theta_0,\theta_1,\dots,\theta_{r-1})$ of partitions with $|\theta_0|+|\theta_1|+\dots+|\theta_{r-1}| = n$}.\ea\] 
We refer to  elements of $\sP_r(n)$ as \emph{$r$-partite partitions of $n$}.  Define  $\psi \wr \lambda$ for $\psi \in \Irr(\ZZ_r)$ and $\lambda \in \sP$  as the character of $\ZZ_r \wr S_{|\lambda|}$ given by
\[ \(\psi \wr \lambda\) (g) = \chi^\lambda(|g|) \cdot (\psi\circ \Delta)(g),\qquad\text{for }g \in \ZZ_r \wr S_{|\lambda|}.\]
 Then, using $\odot$ to denote the external tensor product,  
each irreducible character of $G(r,1,n)$ is of the form 
\[ \chi_\theta \overset{\mathrm{def}}= \Ind_{S_\theta}^{G(r,1,n)} \( \bigodot_{i =0}^{r-1} \psi_i \wr {\theta_i} \),\qquad\text{where }S_\theta = \prod_{i=0}^{r-1} \ZZ_r \wr S_{|\theta_i|},\] for a unique $\theta \in \sP_r(n)$.  The irreducible characters of the wreath product of an arbitrary finite group with $S_n$ arise from a similar construction; see \cite{B91-2} or \cite[Section 4]{M} for a more detailed discussion.

Let $\gamma : G \to \CC$ denote the linear character with 
\be\label{gamma} \gamma(g) = \(\psi_{r/p}\wr (n)\)(g) =  \psi_{r/p} \circ \Delta(g),\qquad\text{for }g\in G.\ee Here $(n)$ denotes the trivial partition of $n$.  Then $\ker \gamma \supset H$ and, since $\gamma$ has order $p$ in the group of linear characters of $G$, it follows that $C = \langle \gamma\rangle = \{ \One,\gamma,\gamma^2,\dots,\gamma^{p-1}\}$ in the notation above.  
A straightforward calculation shows that for all $\theta \in \sP_r(n)$ we have 
\be\label{otimes} \gamma \otimes \chi_\theta = \chi_{\theta'},\qquad\text{where}\qquad \theta'_x = \theta_{x-r/p}\text{ for }x \in [0,r-1],\ee  where with slight abuse of notation we define $\theta_{i-r} = \theta_i$ for $i \in [0,r-1]$. 
If $i \in [0,p-1]$ and $\theta \in \sP_r(n)$, then $\gamma^i\otimes \chi_\theta = \chi_\theta$ if and only if  
\[\theta_j = \theta_{j-ir/p}=\theta_{j-2ir/p} = \dots = \theta_{j-(p-1)ir/p}\] for all $j$.  If this holds then $\frac{p}{\gcd(p,i)}$ is a nontrivial divisor of both $p$ and $n$ since $\sum_{j=0}^{r-1} |\theta_j| = n$.  Hence, if $\gcd(p,n) = 1$, then $\gamma^i \otimes \chi_\theta \neq \chi_\theta$ for all $0<i<p$,  so by the observations above we arrive at the following fact.

\begin{observation}\label{fact} If  $\gcd(p,n)=1$, then each irreducible character of $G(r,p,n)$ is equal to the restriction of exactly $p$ distinct irreducible characters of $G(r,1,n)$.  
\end{observation}

Concerning the irreducible characters of $G(r,p,n)$, we will make use of one additional result due to Caselli.  The next theorem derives from the combination of Proposition 4.4, Theorem 4.5, and Proposition 4.6 in \cite{C2009}.

\begin{theorem}\label{thm4.5} (Caselli \cite{C2009}) Let $r,p,n$ be positive integers with $p$ dividing $r$.  Then 
\[ \left|\left \{ \omega \in G(r,p,n) : \omega^T = \omega \right\} \right | \leq \sum_{\psi \in \Irr\(G(r,p,n)\)} \psi(1)\] and we have equality if and only if $\gcd(p,n) \leq 2$.
\end{theorem}

\subsection{Generalized Involution Models versus Gelfand Models}

The idea of a generalized involution model emerges naturally from the following series of observations.  
Fix a finite group $G$ and for each $\psi \in \Irr(G)$ 
let $\leftexp{\tau} \psi$ denote the irreducible character $\leftexp{\tau}\psi = \psi \circ\tau$.  Define the twisted indicator function $\epsilon_\tau : \Irr(G) \to \{-1,0,1\}$ by 
 \[ \epsilon_\tau(\psi) = \left\{
 \ba 
 1, &\quad \text{if $\psi$ is the character of a representation $\rho$ with $\rho(g) = \overline{\rho(\leftexp{\tau}g)}$ for all $g \in G$,} \\
 0, &\quad \text{if $\psi \neq \overline{\leftexp{\tau}\psi}$}, \\ 
 -1,&\quad\text{otherwise}.\ea\right.\]  When $\tau=1$, this defines the familiar Frobenius-Schur indicator function.  Kawanaka and Matsuyama prove in \cite[Theorem 1.3]{KM1990} that  $\epsilon_\tau$ has the formula
 \[ \epsilon_\tau(\psi) = \frac{1}{|G|} \sum_{g \in G} \psi(g \cdot \leftexp{\tau}{g}),\qquad\text{for }\psi \in\Irr(G).\] 
In addition, we have the following result, which appears in a slightly different form as Theorems 2 and 3 in \cite{BG2004}.  
 
 \begin{theorem} \label{bg-thm} (Bump, Ginzburg \cite{BG2004}) Let $G$ be a finite group with an automorphism $\tau \in \Aut(G)$ 
 such that $\tau^2=1$.  
 Then the following are equivalent:
 \begin{enumerate}
 \item[(i)] The function $\chi : G\to \QQ$ defined by 
  \[ \chi(g) = | \{ u \in G :  u\cdot \leftexp{\tau}{u} = g\}|,\qquad\text{for }g\in G\] is the multiplicity-free sum of all irreducible characters of $G$.
    \item[(ii)] Every irreducible character $\psi$ of $G$ has $\epsilon_\tau(\psi) = 1$.
 \item[(iii)]  The sum $\sum_{\psi \in \Irr(G)} \psi(1)$ is equal to $|\cI_{G,\tau}| = |\{ \omega \in G : \omega \cdot \leftexp{\tau}{\omega} = 1\}|$.
 \end{enumerate}      
\end{theorem}

This theorem motivates Bump and Ginzburg's original definition of a generalized involution model.  In explanation, if the conditions (i)-(iii) hold, then the dimension of any Gelfand model for $G$ is equal to the sum of indices $\sum_i [G : C_{G,\tau}(\omega_i)]$ where $\omega_i$ ranges over a set of representatives of the distinct orbits in $\cI_{G,\tau}$.  The twisted centralizers of a set of orbit representatives in $\cI_{G,\tau}$ thus present an obvious choice for the subgroups $\{H_i\}$ from which to construct a model $\{\lambda_i : H_i \to \CC\}$, and one is naturally tempted to investigate whether $G$ has a generalized involution model with respect to the automorphism $\tau$.

\def\KK{\mathbb{K}}

The following observation concerns the relationship between a generalized involution model and a corresponding \emph{Gelfand model}, which we recall is a representation equivalent to the multiplicity-free sum of all of a group's irreducible representations.  
  Given $\tau \in \Aut(G)$ with $\tau^2=1$ and a fixed subfield $\KK$ of the complex numbers $\CC$, let
\[\cV_{G,\tau} = \KK\spanning\{ C_\omega : \omega \in \cI_{G,\tau}\}\] be a vector space generated by the generalized involutions of $G$.  We often wish to translate a generalized involution model with respect to $\tau\in \Aut(G)$ into a Gelfand model defined in the space $\cV_{G,\tau}$.  
 For this purpose, we repeatedly use the following result, given as Lemma 2.1 in \cite{M}.

\begin{lemma}\label{observation} Let $G$ be a finite group with an automorphism $\tau \in \Aut(G)$ such that $\tau^2=1$.  Suppose 
there exists a function $\sign_G : G\times  \cI_{G,\tau} \to \KK$ such that the map
 $\rho : G\rightarrow \GL(\cV_{G,\tau})$ defined by
\be\label{defn-by} \rho(g) C_\omega = \sign_G(g,\omega)\cdot C_{g\cdot \omega\cdot \leftexp{\tau}g^{-1}},\qquad\text{for }g \in G, \ \omega \in \cI_{G,\tau}\ee  is a representation.  Then the following are equivalent:
\begin{enumerate}
\item[(i)] The representation $\rho$ is a Gelfand model for $G$.
\item[(ii)] The functions 
\[ \left\{ \barr{rccl} \sign_{G}(\cdot, \omega): & C_{G,\tau}(\omega) &\to &\KK  \\ & g & \mapsto & \sign_G(g,\omega) \earr\right\},\] with $\omega$ ranging over any set of orbit representatives of $\cI_{G,\tau}$, form a generalized involution model for $G$.
\end{enumerate}
\end{lemma}

We will apply this lemma to deduce the existence of generalized involution models for $G(r,p,n)$ from certain Gelfand models.  This strategy begins with the following construction for $G(r,1,n)$ due to Adin, Postnikov, and Roichman.  
Fix positive integers $r,n$.  This notation comes from \cite{APR2008}: for any permutation $\pi \in S_n$, define two sets
\[\ba  \Inv(\pi) &= \{ (i,j) : 1\leq i<j\leq n,\ \pi(i) > \pi(j)\}, \\
\Pair(\pi) &= \{ (i,j) : 1\leq i < j \leq n,\ \pi(i) = j,\ \pi(j) = i\}.\ea\]  The set $\Inv(\pi)$ is the inversion set of $\pi$ and its cardinality is equal to the element's length, by which we mean the minimum number of factors needed to write $\pi$ as a product of simple reflections.  In particular, the signature of $\pi$ is $\sgn(\pi) =(-1)^{|\Inv(\pi)|}$.  The set $\Pair(\pi)$ corresponds to the set of cycles of length two in $\pi$.
   Next, for any two elements $g,\omega \in G(r,1,n)$, let $B(g,\omega)$ denote the subset of $[1,n]$ given by 
\[ B(g,\omega) =\left\{\barr{ll} \varnothing, &\text{if $r$ is odd}, \\  \\
\left\{ i \in \Fix(|\omega|) 
: \ba 
&z_\omega(i) \text{ is odd and } z_g(i)+k\in[r/2,r-1] \\
& \text{for the $k \in [0,r/2-1]$ with $2k+1=z_\omega(i)$}\ea
\right\}, &\text{if $r$ is even.}\earr\right.\]  Here $\Fix(\pi) =\{ i\in [1,n] : \pi(i) = i\}$ denotes the set of fixed points of a permutation $\pi \in S_n$.  

Now let $\cV_{r,n}$ be a vector space generated by the symmetric elements in $G(r,1,n)$,
\[\cV_{r,n} = \QQ\spanning\left\{ C_\omega : \omega \in G(r,1,n),\ \omega^T = \omega\right\},\] 
and define $\rho_{r,n} :G(r,1,n) \to \GL(\cV_{r,n})$ by
\[ \rho_{r,n}(g) C_\omega ={\sign_{r,n}(g, \omega)}_{} \cdot C_{g \omega g^T},\qquad\text{for }g,\omega \in G(r,1,n) \text{ with }\omega^T = \omega.\] 
where
\[ \sign_{r,n}(g,\omega) = (-1)^{|B(g,\omega)|}\cdot (-1)^{|\Inv(|g|) \cap \Pair(|\omega|)|},\qquad\text{for }g,\omega \in G(r,1,n).\]  
This map defines a representation; in fact, we have the following result given as Theorem 1.2 in \cite{APR2008}.

\begin{theorem} \label{thm1.2} (Adin, Postnikov, Roichman \cite{APR2008}) The map $\rho_{r,n}$ is a Gelfand model for $G(r,1,n)$.
\end{theorem}
Observe that $\cV_{r,n} = \cV_{G, \tau}$ where $\KK = \QQ$ and $G=G(r,1,n)$ and $\tau$ is the inverse transpose automorphism  $\tau : g\mapsto \overline g $.  Thus, it follows from Lemma \ref{observation} that $G(r,1,n)$ has a generalized involution model with respect to $\tau$.  Theorems 5.2 and 5.3 in \cite{M} give an explicit description of this model and discuss some of its properties.

In the next section we modify the representation $\rho_{r,n}$ to define Gelfand models for $G(r,p,n)$ when $\gcd(p,n) = 1$.  For this purpose, we state now a needed lemma.
Assume $r$ is even.  We then have two $\rho_{r,n}$-invariant subspaces of $\cV_{r,n}$ given by
\be\label{+-}\ba  \cV_{r,n}^+ &\overset{\mathrm{def}}= \QQ\spanning\left\{ C_\omega : \omega \in G(r,1,n),\ \omega^T = \omega,\ \Delta(\omega) \in 2\ZZ_r \right\}, \\
 \cV_{r,n}^- &\overset{\mathrm{def}}= \QQ\spanning\left\{ C_\omega : \omega \in G(r,1,n),\ \omega^T = \omega,\ \Delta(\omega) \notin 2\ZZ_r \right\}. \ea\ee 
 Let $\chi_{r,n}^+$ and $\chi_{r,n}^-$ denote the characters of $G(r,1,n)$ corresponding to the subrepresentations of $\rho_{r,n}$ on $\cV_{r,n}^+$ and $\cV_{r,n}^-$ respectively.  Recall from (\ref{gamma}) above the definition of the linear character $\gamma$ of $G(r,1,n)$. The following lemma, given as Proposition 5.2 in \cite{M}, will play a useful part in our proofs in Section \ref{positive}.
  
\begin{lemma} \label{plusminus-cor} Let $r,p,n$ be positive integers with $r$ even and $p$ dividing $r$. Then
\[ \gamma\otimes \chi_{r,n}^+= \left\{\ba & \chi_{r,n}^-, &&\text{if $n$ and $r/p$ are odd}, \\
&\chi_{r,n}^+,&&\text{otherwise},\ea\right.
\qquad
\gamma\otimes \chi_{r,n}^- =\left\{\ba & \chi_{r,n}^+, &&\text{if $n$ and $r/p$ are odd}, \\
&\chi_{r,n}^-,&&\text{otherwise}.\ea\right.
\]

\end{lemma}

\section{Constructions for the Infinite Series}\label{positive}

In this section we describe how one can obtain a generalized involution model for  $G(r,p,n)$ in the cases where this is possible.  We have two methods for doing this: by explicitly identifying the set of linear characters comprising our model, or by giving a Gelfand model of the particular form appearing in Lemma \ref{observation}.  We apply the second method when $\gcd(p,n) =1$ and the first when $\gcd(p,n)=n=2$ and $r/p$ is odd. In Section {\ref{classifying-sec} we will discover that $G(r,p,n)$ does not have a generalized involution model in any other cases.

\subsection{Gelfand Models for $G(r,p,n)$ with $p,n$ Coprime}

Some work has been done on this topic.
In the recent paper \cite{C2009}, Caselli describes a representation for $G(r,p,n)$ in  the complex vector space 
\[ \cV_{r,p,n}^\CC \overset{\mathrm{def}} = \CC\spanning\left\{ C_\omega : \omega \in G(r,p,n),\ \omega^T = \omega\right\}\] which defines a Gelfand model whenever equality obtains in Theorem \ref{thm4.5}; i.e., when $\gcd(p,n) \leq 2$.  He calls such complex reflection groups \emph{involutory}.  Caselli's constructions 
do not arise from generalized involution models and require the field of complex numbers for their definition.  The Gelfand models we present coincide with Caselli's only when $r=1$, and are by contrast rational representations.  

To give these, we begin by noting that the Gelfand model $\rho_{r,n}$ for $G(r,1,n) $ restricts to a representation of $G(r,p,n)$ for any $p$ dividing $r$, and that one obvious subrepresentation of this restriction poses a natural candidate for a Gelfand model.  Specifically, if we define $\cV_{r,p,n} \subset \cV_{r,n}$ as the subspace
\[ \cV_{r,p,n} =  \QQ\spanning\left\{ C_\omega : \omega \in G(r,p,n),\ \omega^T = \omega\right\}\] then since $G(r,p,n)$ is closed under taking transposes$-$as defined by (\ref{book})$-$the map
$\rho_{r,p,n} : G(r,p,n) \to \GL(\cV_{r,p,n})$ given by
\[ \rho_{r,p,n}(g) C_\omega = \sign_{r,n}(g,\omega) \cdot C_{g\omega g^T},\qquad\text{for }g,\omega \in G(r,p,n) \text{ with }\omega^T = \omega\] is automatically a well-defined $G(r,p,n)$-representation.  
The following theorem says exactly when this representation is a Gelfand model.

\begin{theorem}\label{reflection-thm} Let $r,p,n$ be positive integers with $p$ dividing $r$.  Then the representation $\rho_{r,p,n}$ is a Gelfand model for $G(r,p,n)$ if and only if $\gcd(p,n) = 1$ and $p$ or $r/p$ is odd.  
\end{theorem}

\def\cU{\mathcal{U}}
\def\cW{\mathcal{W}}

\begin{proof}
Let $G= G(r,1,n)$ and $H = G(r,p,n)$.  
By Theorem \ref{thm4.5}, $\rho_{r,p,n}$ can only be a Gelfand model for $H$ if $\gcd(p,n) = 1$ or 2, so we only consider those cases.
View $\cV_{r,n}$ as a $G$-module by defining $g C_\omega  = \rho_{r,n}(g) C_\omega$ for $g \in G$, and for any $i \in \ZZ_r$, let $\cV_{r,n}(i)$ denote the $H$-submodule
\[\cV_{r,n}(i) = \QQ\spanning\left\{ C_\omega : \omega \in G,\ \omega^T = \omega,\ \Delta(\omega) - i\in p\ZZ_r\right\}.\] 
Observe that $\cV_{r,p,n}= \cV_{r,n}(0)$ and that $\cV_{r,n} 
= \cV_{r,n}(0) \oplus \cV_{r,n}(1)\oplus \cdots \oplus \cV_{r,n}(p-1)$.
 
  Suppose $\gcd(p,n) = 1$, and  let $c  \in G$  denote the central element 
  \be\label{c_i}c  = \( (1,1,\dots,1), 1\) \in G\qquad\text{so that}\qquad c^i = \( (i,i,\dots,i),1\).\ee   Observe that $\cV_{r,n}(2ni+j) = c^i  \cV_{r,n}(j)$ and so $\cV_{r,n}(j) \cong \cV_{r,n}(2ni+j)$ as $H$-modules since $c^i$ is central.  Consequently, if $p$ is odd then $\gcd(p,2n) =1$ and the $H$-modules  $\cV_{r,n}(i)$ are all isomorphic, since as $i$ ranges over $0,1,\dots,p-1 \in \ZZ_r$, the elements $2ni$ represent every coset of $p\ZZ_r$ in $\ZZ_r$.  In this case, it follows that an irreducible $H$-module $\cU$ is a constituent of $\cV_{r,p,n} = \cV_{r,n}(0)$ with multiplicity $m$ if and only if $\cU$ is a constituent of $\cV_{r,n}$ with multiplicity $ pm$.  Therefore if $p$ is odd then $\rho_{r,p,n}$ is a Gelfand model for $H$ since $\gcd(p,n)=1$ implies that each irreducible $H$-module appears as a constituent of $\cV_{r,n}$ with multiplicity $p$ by Observation \ref{fact}.
   
  Suppose alternatively that $\gcd(p,n) = 1$ but $p$ is even, so that $n$ is odd.  Then 
 by the same considerations the $H$-modules $\cV_{r,n}^+$ and $\cV_{r,n}^-$ defined by (\ref{+-}) 
are isomorphic to $p/2$ copies of $\cV_{r,n}(0)$ and $\cV_{r,n}(1)$, respectively. 
Since every irreducible $H$-module is isomorphic to a constituent of $\cV_{r,n}$ with multiplicity $p$ as $\gcd(p,n) = 1$, it follows that every irreducible $H$-module is isomorphic to a constituent of $\cV_{r,p,n} = \cV_{r,n}(0)$ with multiplicity one if and only if $\cV_{r,n}^+ \cong \cV_{r,n}^-$ as $H$-modules.  
We claim that this holds if and only if $r/p$ is odd.  

To show this, observe that $\cV_{r,n}^+ \cong \cV_{r,n}^-$ as $H$-modules if and only if $\gamma\otimes \chi_{r,n}^+ = \chi_{r,n}^-$, where $\gamma$ is the character defined by (\ref{gamma}).  The ``if'' direction of this statement is immediate since $\gamma$ restricts to the trivial character of $H$, and the other direction follows from Lemma \ref{plusminus-cor}, since if $\gamma\otimes \chi_{r,n}^+ \neq \chi_{r,n}^-$ then $\gamma \otimes \chi_{r,n}^- = \chi_{r,n}^-$ which implies that no irreducible constituent of the nonzero $H$-module $\cV_{r,n}^-$ appears as a constituent of $\cV_{r,n}^+$.  Since $n$ is odd, Lemma \ref{plusminus-cor} implies that $\gamma\otimes \chi_{r,n}^+ = \chi_{r,n}^-$ if and only if $r/p$ is odd, which proves our claim.   Thus if $\gcd(p,n) =1$ then $\rho_{r,p,n}$ is a Gelfand model for $H$ if and only if $p$ or $r/p$ is odd.

To complete the proof, suppose $\gcd(p,n)=2$ so that $n$ and $p$ are both even.  Then $r$ is even and it follows from Lemma \ref{plusminus-cor} that $\gamma\otimes \chi_{r,n}^- = \chi_{r,n}^-$.  Hence any irreducible constituent of the nonzero $H$-module $\cV_{r,n}^-$ does not appear as a constituent of $\cV_{r,n}^+$, or in the submodule $\cV_{r,p,n}  = \cV_{r,n}(0)$, so $\rho_{r,p,n}$ cannot be a Gelfand model for $H$.
\end{proof}

Suppose $\gcd(p,n)=1$ but both $p$ and $r/p$ are even.  Then while Theorem \ref{reflection-thm} does not hold, by modifying our construction slightly we can still produce a Gelfand model for $G(r,p,n)$ in $\cV_{r,p,n}$.  In this case, for every $\omega \in G(r,p,n)$ exactly one of the containments $\Delta(\omega) \in 2p\ZZ_r$ or $\Delta(\omega) - p \in 2p \ZZ_r$ holds.  Thus, we may define $\wt B(g,\omega)$  for two elements $g,\omega \in G(r,p,n)$ as the subset of $[1,n]$ given by 
\[ \wt B(g,\omega) =\left\{\barr{ll} 
\left\{ i \in \Fix(|\omega|) 
: \ba 
&
z_\omega(i) \text{ is odd and } z_g(i)+k\in[r/2,r-1] \\
&
 \text{for the $k \in [0,r/2-1]$ with $2k+1=z_\omega(i)$}
 \ea
\right\}, &\text{if $\Delta(\omega) \in 2p\ZZ_r$},\\
\\
\left\{ i \in \Fix(|\omega|) 
: \ba  
&
z_\omega(i) \text{ is even and } z_g(i)+k\in[r/2,r-1] \\
&
 \text{for the $k \in [0,r/2-1]$ with $2k=z_\omega(i)$}
 \ea
\right\}, &\text{if $\Delta(\omega) \notin 2p\ZZ_r$}.
\earr\right.\]  
Let $\wt \rho_{r,p,n} : G(r,p,n) \to \GL(\cV_{r,p,n})$ be the map given by 
\[ \wt \rho_{r,p,n}(g) C_\omega = \wt\sign_{r,p,n}(g,\omega) \cdot C_{g\omega g^T},\qquad\text{for $g,\omega \in G(r,p,n)$ with $\omega^T = \omega$}\] where 
\[ \wt \sign_{r,p,n}(g,\omega) =  (-1)^{|\wt B(g,\omega)|}\cdot (-1)^{|\Inv(|g|) \cap \Pair(|\omega|)|},\qquad\text{for }g,\omega \in G(r,1,n).\] 
We then have the following result.  
\begin{theorem}\label{reflection-thm2}  
Let $r,p,n$ be positive integers with $p$ dividing $r$.  If $\gcd(p,n) =1$ but $r/p$ and $p$ are both even,   then $\wt \rho_{r,p,n}$ is a Gelfand model for $G(r,p,n)$.
\end{theorem}

In the following proof, it is helpful to note that if $c \in G(r,1,n)$ is the central element defined by (\ref{c_i}), then 
\[ \wt \sign_{r,p,n}(g,\omega) = \left\{\ba& \sign_{r,n}(g,\omega)\cdot C_{g \omega g^T}, && \text{if }\Delta(\omega) \in 2p\ZZ_r, \\
&\sign_{r,n}(g,\omega c)\cdot C_{g \omega g^T}, && \text{if }\Delta(\omega) -p \in 2p\ZZ_r. \ea\right.\]
  Given this observation, one can check without difficulty that $\wt \rho_{r,p,n}$ is a well-defined representation when $p$ and $r/p$ are both even, using the fact that $\rho_{r,p,n}$ is a representation and that $\Delta(g\omega g^T) - \Delta(\omega) = 2\Delta(g) \in 2p\ZZ_r$ for all $g,\omega \in G(r,p,n)$.

\begin{proof}
Again let $G= G(r,1,n)$ and $H = G(r,p,n)$, and view $\cV_{r,p}$ as a $G$-module as in the proof of Theorem \ref{reflection-thm}.  Since $2p$ divides $r$, $\cV_{r,p}$ decomposes into a direct sum of $2p$ distinct $H$-submodules as $\cV_{r,n} = \wt \cV_{r,n}(0) \oplus \wt \cV_{r,n}(1)\oplus \dots \oplus \wt \cV_{r,n}(2p-1)$ where \[\wt\cV_{r,n}(i) = \QQ\spanning\left\{ C_\omega : \omega \in G,\ \omega^T = \omega,\ \Delta(\omega)- i \in 2p \ZZ_r\right\}.\]
Defining $c \in G$ by (\ref{c_i}), we again have $\wt \cV_{r,n}(2ni+j) = c^i \wt \cV_{r,n}(j)$ for all $i,j \in \ZZ_r$, so since $\gcd(2p,n) = 1$ as $p$ is even, the $H$-modules $\cV_{r,n}^+$ and $\cV_{r,n}^-$ defined by (\ref{+-}) 
are isomorphic to $p$ copies of $\wt \cV_{r,n}(0)$ and $\wt \cV_{r,n}(1)$, respectively.    Since $\delta\otimes \cV_{r,n}^\pm = \cV_{r,n}^\pm$ by Lemma \ref{plusminus-cor}, the $H$-modules $\cV_{r,n}^+$ and $\cV_{r,n}^-$ do not share any irreducible constituents.  Therefore, since each irreducible $H$-module appears as a constituent of $\cV_{r,n}$ with multiplicity $p$ by Observation \ref{fact},  it follows that each irreducible $H$-module appears as a constituent of $\wt \cV_{r,n}(0) \oplus \wt \cV_{r,n}(1)$ with multiplicity one.

If we view $\cV_{r,p,n}$ as an $H$-module by defining $g C_\omega = \wt \rho_{r,p,n}(g) C_\omega$ for $g \in H$, then $\cV_{r,p,n}$ decomposes into $H$-submodules as $\cV_{r,p,n}  = \cV_0 \oplus \cV_1$ where 
\[ \cV_0 = \QQ\spanning \{ C_\omega : \Delta(\omega) \in 2p\ZZ_r \}\qquad\text{and}\qquad
\cV_1 = \QQ\spanning \{ C_\omega : \Delta(\omega) -p\in 2p\ZZ_r\}.\]
By definition $\cV_0= \wt \cV_{r,n}(0)$, and 
one easily sees that the linear map $\cV_1 \to \wt \cV_{r,n}(p+n)$ defined on basis elements by $C_\omega \mapsto C_{\omega c}$  is an isomorphism of $H$-modules.  Since $n$ is odd, $\wt\cV_{r,n}(p+n) \cong \wt \cV_{r,n}(1)$ as $H$-modules, and so we conclude that $\wt \rho_{r,p,n}$ is a Gelfand model.
\end{proof}

Since in the notation of the previous section $\cV_{r,p,n}$ is precisely the vector space $\cV_{G,\tau}$ with $\KK = \QQ$, $G=G(r,p,n)$, and $\tau \in \Aut(G)$ the inverse transpose automorphism $\tau : g \mapsto \overline g$, 
 we are  afforded the following corollary by Lemma \ref{observation}.

\begin{corollary}\label{thm-cor}
Let $r,p,n$ be positive integers such that $p$ divides $r$.  Then $G(r,p,n)$ has a generalized involution model  with respect to the inverse transpose automorphism $g\mapsto \overline g$ if $\gcd(p,n)=1$.
\end{corollary}

\begin{remark}
One can form a generalized involution model for $G=G(r,p,n)$ by choosing a set of representatives $\{\omega\}$ for the $\tau$-twisted conjugacy classes in $\cI_{G,\tau}$, and then taking the linear characters $\lambda : C_{G,\tau}(\omega) \to \QQ$ defined as the coefficients in $\QQ$ such that $\rho(g) C_\omega = \lambda(g) C_\omega$ for those $g \in G$ with $g\omega g^T = \omega$, where $\rho$ is our Gelfand model.  
\end{remark}

\subsection{A Generalized Involution Model for $G(r,p,2)$}

We can only expect to be able to construct a generalized involution model for $G(r,p,n)$ when $\gcd(p,n) \leq 2$, and we will in fact be unable to do so  when $\gcd(p,n)=2$ in most cases.  Here we deal with the one exception to this rule, occurring when $n=2$ and $r/p$ is odd.  In contrast to the previous section, here we produce the generalized involution model directly.

Throughout this section, fix positive even integers $r,p$ with $p$ dividing $r$ such that $r/p$ is odd. We write $G = G(r,p,2)$ and let $\tau \in \Aut(G)$ denote the inverse transpose automorphism  $\tau : g\mapsto \overline g$.  Of immediate relevance is the following consequence of Theorem \ref{bg-thm}.

\begin{lemma} \label{model-char}Let $g = \( (a,b), \pi \) \in G$, so that  $a,b \in \ZZ_r$ such that $a+b \in p\ZZ_r$ and $\pi \in S_2$.  Then 
\be\label{model-char-form} \sum_{\psi \in \Irr\(G\)} \psi(g) = \left\{\barr{ll} ({r^2+2 r})/{p}, &\text{if }a=b=0\text{ and }\pi=1, \\ {2 r}/{p},&\text{if }a=-b \in 2\ZZ_r \setminus\{0\}\text{ and }\pi = 1, \\ 0,&\text{otherwise}.\earr\right.\ee
\end{lemma}

\begin{proof}
By Theorems \ref{thm4.5} and \ref{bg-thm}, it suffices to show that the right-hand side is equal to the number of $\omega \in G$ with $\omega\cdot\leftexp{\tau}\omega = g$.  Let $\omega = \( (x,y), \sigma \) \in G$; then $(x,y) \in (\ZZ_r)^2$ can assume $r^2/p$ distinct values.     If $\sigma = 1$ then $\omega \cdot \leftexp{\tau}\omega = \((0,0),1\)$ while if $\sigma\neq 1$ then $\omega \cdot \leftexp{\tau}\omega = \( (y-x,x-y),1\)$.  As there are $2 r/p$ choices of $(x,y) \in (\ZZ_r)^2$ such that $x+y \in p \ZZ_r$ and $x-y = b$ if $a,b \in 2\ZZ_r$ and zero choices if $a,b \notin 2\ZZ_r$, the lemma follows. 
\end{proof}

Let $s_1 \in S_2$ denote the simple reflection $s_1 = (1\cyc 2)$.  One checks that the elements
\[\omega_1 = \( (0,0), 1\),\qquad \omega_2 = \( (1,-1), 1\),\qquad \omega_3 = \( (0,0), s_1\)\qquad \omega_4 = \( (p/2, p/2), s_1\)\] represent the distinct $\tau$-twisted conjugacy classes in $\cI_{G,\tau}$, and that 
\[\ba &\barr{c} C_{G,\tau}(\omega_1) = \left\{ \( 0,1\), \((\frac{r}{2},\frac{r}{2}),1\), \( 0, s_1\), \((\frac{r}{2},\frac{r}{2}),s_1\)\right\} \cong S_2 \times S_2,
\earr\\
&\barr{c} C_{G,\tau}(\omega_2) = \left\{ \( 0,1\),  \((\frac{r}{2},\frac{r}{2}),1\), \( (-1,1), s_1\),\((\frac{r}{2}-1,\frac{r}{2}+1),s_1\)\right\} \cong S_2 \times S_2,\earr
\ea\]
and $C_{G,\tau}(\omega_3) = C_{G,\tau}(\omega_4) = G(r,r,2).$   Define linear characters $\lambda_i : C_{G,\tau}(\omega_i) \to \QQ$ by 
\[ \barr{c|cccc} 
 &\( (0,0),1\)  & \((\frac{r}{2},\frac{r}{2}),1\) & \( (-1,1), s_1\) & \((\frac{r}{2}-1,\frac{r}{2}+1),s_1\) \\
 \hline \lambda_2 & 1 & -1 & -1 & 1
\earr\] and 
\[
\lambda_1(g) =1,\qquad \lambda_3(g) = \sgn(|g|),\qquad \lambda_4(g) = \sgn(|g|) \cdot (-1)^{z_g(1)}.\] 
In the definition of $\lambda_4$ we are of course viewing $z_g(1) \in \ZZ_r$ as an integer in $[0,r-1]$; the given formula only makes sense because $n=2$.

We now have the following result.

\begin{proposition}
Let $r,p$ be even positive integers with $p$ dividing $r$, such that $r/p$ is odd.  Then the linear characters $\lambda_i : C_{G,\tau}(\omega_i) \to \QQ$ for $1\leq i\leq 4$ form a generalized involution model for $G = G(r,p,2)$.
\end{proposition}

\begin{proof}
If we define $h_{ij} = \( ( ip+j,-j),1\) \in G$ for $i,j \in \ZZ_r$ and let
\[\ba
\cC_1 &= \cC_2 = \{h_{ij} : i \in [0,r/p-1],\ j \in [0,r/2-1]\}, &&\qquad \text{so that }|\cC_1|=|\cC_2| = r^2/(2p)\\
 \cC_3 &= \cC_4 = \{ h_{i0} : i \in [0,r/p-1]\},&&\qquad \text{so that } |\cC_3|=|\cC_4| = r/p
 \ea\]  then each $\cC_i$ forms a set of left coset representatives of $C_{G,\tau}(\omega_i)$ in $G$.   Let $g = \( (a,b), \pi \) \in G$ denote an arbitrary element of $G$ with $a,b \in \ZZ_r$, $\pi \in S_2$ and $a+b \in p\ZZ_r$.  Observe that 
 \[ h_{ij} \cdot g\cdot  (h_{ij})^{-1} = \left\{\barr{ll} \( (a,b),1\), & \text{if }\pi = 1, \\ 
 \( ( a-ip-2j, b+ip+2j), s_1\),&\text{if }\pi \neq 1.\earr\right.\] Write $\Lambda_i = \Ind_{C_{G,\tau}(\omega_i)}^G(\lambda_i)$.  Using the preceding observation with the Frobenius formula for induced characters, it is not difficult to check that:
 \begin{enumerate}
 \item[(i)] If $a+b\neq 0$ then $\Lambda_i (g) = 0$ for $1\leq i\leq 4$.
 
 \item[(ii)] If $a+b =0$ and $\pi=1$ then 
 \[\barr{lcl}
  \Lambda_1(g)=
 \left\{ \barr{ll} r^2/(2p), & \text{if }a=b \in \{0,r/2\}, \\ 
 0,&\text{otherwise,}\earr
 \right. 
 & &
 \Lambda_3(g) = r/p, 
 \\ \\
 \Lambda_2(g) = 
 \left\{ \barr{ll} r^2/(2p), & \text{if }a=b =0, \\ 
 - r^2/(2p), & \text{if }a=b =r/2, \\ 
 0,&\text{otherwise,}\earr
 \right.
 & &
  \Lambda_4(g) = 
 \left\{ \barr{ll} r/p, & \text{if }a \in 2\ZZ_r, \\ 
 - r/p, & \text{if }a \notin 2\ZZ_r.
 \earr\right.
 \earr
\]

 \item[(iii)] If $a+b =0$ and $\pi\neq 1$ then 
 \[\barr{lcl}
  \Lambda_1(g)=
 \left\{ \barr{ll} 2r/p, & \text{if }a \in 2\ZZ_r\text{ and } r/2 \text{ is even}, \\ 
  0, & \text{if }a \notin 2\ZZ_r\text{ and } r/2 \text{ is even}, \\ 
 r/p, & \text{if } r/2 \text{ is odd}, 
 \earr
 \right. 
 & &
 \Lambda_3(g) = -r/p, 
 \\ \\
 \Lambda_2(g) = 
 \left\{ \barr{ll} 0, & \text{if } r/2 \text{ is even}, \\ 
 r/p, & \text{if }a \in 2\ZZ_r\text{ and } r/2 \text{ is odd}, \\ 
 -r/p, & \text{if }a \notin 2\ZZ_r\text{ and } r/2 \text{ is odd}, 
 \earr
  \right.
 & &
  \Lambda_4(g) = 
 \left\{ \barr{ll} -r/p, & \text{if }a \in 2\ZZ_r, \\ 
 r/p, & \text{if }a \notin 2\ZZ_r.
 \earr\right.
 \earr
\]   
 
 \end{enumerate}
In turn, these formulas imply that $\(\Lambda_1 + \Lambda_2 + \Lambda_3 + \Lambda_4\)(g)$ is precisely equal to the right-hand side of (\ref{model-char-form}), which completes our proof.
\end{proof}

\section{Automorphisms of $G(r,p,n)$}

To prove that the groups $G(r,p,n)$ do not have generalized involution models other than in the situations addressed by the previous section, we require some understanding of these groups' automorphisms.  
In particular, we require a sufficiently explicit description of the elements of $\Aut\(G(r,p,n)\)$ to be able to deduce precisely which automorphisms can satisfy the conditions of Theorem \ref{bg-thm}.

Marin and Michel  provide in \cite{MM} several useful general results concerning the structure of $\Aut(G)$ when $G$ is any finite complex reflection group.  In particular, they prove that when $G$ is an irreducible complex reflection group not equal to the symmetric group $S_6$, each automorphism of $G$ is the composition of an automorphism which preserves the pseudo-reflections in $G$ and a ``central automorphism,'' by which we mean a map $\tau$ such that $\leftexp{\tau} g\cdot g^{-1}$ is always central.  Letting $V$ denote the vector space on which $G$ acts irreducibly, Marin and Michel describe how  each reflection-preserving automorphism can be interpreted as the composition of an automorphism induced from the normalizer of $G$ in $\GL(V)$ and an automorphism induced from the Galois group of $\KK$ over $\QQ$, where $\KK$ is the field of definition, i.e., the extension of $\QQ$ generated by the traces of elements of $G$.  They further discuss how to construct central automorphisms from the linear characters of $G$.

Marin and Michel's paper does not go as far as to actually write down the definitions of all the automorphisms in a very accessible fashion.
 Shi and Wang  in the article \cite{ShiWang} do write down explicit formulas, but only for the subgroup of reflection-preserving automorphisms of $G(r,p,n)$.
From elementary considerations and without too much difficulty, one can give a complete and explicit description of $\Aut\(G(r,p,n)\)$, and we provide this here for completeness. It is possible to glean many of these results from Shi and Wang's classification \cite{ShiWang} and Marin and Michel's work \cite{MM}.  The content of this section is, as such, to produce from a short, self-contained argument actual formulas for all the automorphisms of $G(r,p,n)$ which we can use to classify their generalized involution models.

In what follows, we denote by $\Inn(G)$  the group of inner automorphisms of a group $G$; by $\Out(G)$ the quotient group $\Aut(G) / \Inn(G)$; and by $Z(G)$ the center of $G$.
Fix positive integers $r,p,n$ with $p$ dividing $r$.  
Let $e_i = (0,\dots,0,1,0,\dots,0) \in (\ZZ_r)^n$ denote the standard vector in the obvious free basis of $(\ZZ_r)^n$ over $\ZZ_r$, and define  elements $s_i,s_i',s\in G(r,r,n)$ and $t,c \in G(r,1,n)$ by
\be\label{gen} \ba &s_i =\(0,  (i\cyc i+1)\)  \text{ for $i=1,\dots,n-1$,} \\
&s_i'  =\( e_i-e_{i+1}, (i\cyc i+1)\)\text{ for $i=1,\dots,n-1$,} \\ 
&s= \( e_1-e_2, 1\),  \\
&t = \(e_1,1\), \\
&c =\( e_1 +  e_2 + \dots + e_n,1\) .\ea\ee  Note that the elements $s_i'$, $s$, are only defined for $r\geq 2$ and $n\geq 2$, and that $s_1' = s_1s$.  Also, observe that each $s_i$ and $s_i'$ has order 2 while $s$, $t$, $c$ all have order $r$.  In particular when $r=1$ we have $s=t=c=1$. 

The group $G(r,p,n)$ is generated by $s_1'$, $s_1$, $\dots$, $s_{n-1}$, $t^p$ or by $s_1$, $\dots$, $s_{n-1}$, $s$, $t^p$; we can omit from these lists  $s_1'$ and $s$ if $p=1$ and $t^p$ if $p=r$.  Brou\'e, Malle, and Rouquier \cite{BMR} give presentations for $G(r,p,n)$, as well as for the exceptional groups $G_i$, using the former set of generators; however, the latter set will be more convenient in many of our definitions.  

For each integer $j$ we note that
\be\label{c-eqn} c^{j} = t^j\cdot \( s_1 t ^js_1\) \cdot \(s_2s_1 t^j s_1 s_2\) \cdots \( s_{n-1}\cdots s_2 s_1 t^j s_1 s_2 \cdots s_{n-1}\) \in Z\( G(r,p,n)\).\ee  
The center of $G(r,p,n)$ almost always lies in the subgroup of $G(r,1,n)$ generated by $c$, as we note for later use in the following basic lemma.

\begin{lemma}\label{center}
Let $r,p,n$ be integers with $p$ dividing $r$.  
 If  $d = \gcd(p,n)$ then 
\[ (\ZZ_r)^n \cap Z\(G(r,p,n)\) = \{ c^{jp/d} : j \in [0,dr/p-1]\}.\]  This subgroup is equal to the center of $G(r,p,n)$ unless $(r,p,n)$ is $(1,1,2)$ or $(2,2,2)$, in which case $G(r,p,n)$ is abelian.

\end{lemma}

\begin{proof}
We leave this easy exercise to the reader.
%
%
\end{proof}

To make our notation less cumbersome, we set 
\[ C(r,p,n) = (\ZZ_r)^n \cap Z\(G(r,p,n)\).\]  The following definition in some sense names all nontrivial outer automorphisms of $G(r,p,n)$.  
Given $j,k \in \ZZ$  and $z \in C(r,1,n)$, let $\alpha_{j,k,z} : G(r,1,n) \to G(r,1,n)$ be the map 
\be \label{alpha} \alpha_{j,k,z} : (x,\pi) \mapsto z^{\ell(\pi)} c^{\Delta(x)\cdot k} (jx,\pi),\qquad\text{for }x \in (\ZZ_r)^n,\ \pi \in S_n.\ee  We recall that $\Delta : (\ZZ_r)^n \to \ZZ_r$ is the homomorphism $\Delta(x) = x_1 + x_2 + \dots x_n$.  In our superscripts we
 naturally identify $\ZZ_r$ with the integers $[0,r-1]$ and view $\ZZ_r$ as a $\ZZ$-module.  Since $c$ has order $r$, this is well defined.  Also,  $\ell : S_n \to \ZZ_{\geq 0}$ denotes the usual length function, defined as the minimum number of factors needed to write a permutation as a product of the simple transpositions $s_i$, or equivalently the cardinality of a permutation's inversion set.  
 
 The map $\alpha_{j,k,z}$ has the following effect on our generators: 
 \[ s_i  \mapsto  zs_i, \qquad s_i' \mapsto z s_i'(s_is_i')^{j-1},\qquad
  s \mapsto  s^j, \qquad
 t \mapsto  c^{k} t^j, \qquad
  c \mapsto  c^{j+nk}.\]  Observe that $\alpha_{1,0,1}$ is the identity and $\alpha_{-1,0,1}$ is the inverse transpose automorphism $g\mapsto \overline g$. 
    The map $\alpha_{j,k,z}$ is often but not always an automorphism, as we see in the following lemma.

\begin{lemma}\label{aut-prop} Let $r,p,n$ be positive integers with $p$ dividing $r$.   If $j,k \in \ZZ$ and $z \in C(r,1,n)$, then the map $\alpha_{j,k,z}$ restricts to an automorphism of $G(r,p,n)$ if and only if 
\be\label{aut-prop-cond} \gcd(j,r) = \gcd(j+nk,r/p) = 1\qquad\text{and} \qquad z \in C(r,p,n)\text{ and }z^2=1.\ee
\end{lemma}

\begin{proof}
Assume $\alpha_{j,k,z}$ restricts to an automorphism of $G(r,p,n)$.   The image $z s_i$ of $s_i$ then has order two and belongs to $G(r,p,n)$, so $z^2=1$ and $z = zs_i \cdot s_i \in G(r,p,n)$, which implies  
 $z \in C(r,p,n)$  since $z \in C(r,1,n)$. Likewise, the image $s^j$ of $s$ has order $r$ so $\gcd(j,r)=1$, and the image  $c^{p(j+nk)}$ of $c^p$
has order $r/p$ so $\gcd(j+nk,r/p)=1$.

Conversely, suppose (\ref{aut-prop-cond}) holds. 
Since $c^{p} \in G(r,p,n)$ and $\Delta(jx) = j \cdot \Delta(x)$,  it follows that $\alpha_{j,k,z}$ maps $G(r,p,n)$ into itself.  
One easily checks that $\alpha_{j,k,z}$ is a homomorphism using the following observations:  $c,z$ are central; if $\pi,\sigma \in S_n$ then $\ell(\pi) + \ell(\sigma) - \ell(\pi\sigma)$ is even; and $\Delta( \sigma(x)) = \Delta(x)$ for all $x \in (\ZZ_r)^n$ and $\sigma \in S_n$.  It remains only to show that $\alpha_{j,k,z} : G(r,p,n) \to G(r,p,n)$ is bijective, and for this it suffices to show that $\alpha_{j,k,z}$ is surjective.

To prove this, we first observe that $z$ is either the identity or the element $c^{r/2}$ when $r$ is even and $nr/2$ is a multiple of $p$.  Assume this latter case occurs; since  $\Delta\(\frac{r}{2} e_1 + \dots + \frac{r}{2}e_n\) = nr/2$ we then have $\alpha_{j,k,z}( z)=  c^{(j+nk)r/2}.$     If $r/p$ is even then $\gcd(j+nk,r/p)=1$ implies that  $j+nk$ is odd.  If $r/p$ is odd then $p$ is even and $r/2$ is an odd multiple of $p/2$, so $n$ must be even in order for $nr/2$ to be a multiple of $p$.  Since $r$ is even and $\gcd(j,r)=1$, $j$ is odd, so again $j+nk$ is odd.  Hence 
 $c^{(j+nk)r/2} = c^{r/2} = z$, and we conclude that in either case $\alpha_{j,k,z}(z) = z$.  
 
Given this observation, it follows that $\alpha_{j,k,z} (z s_i)= s_i$ for all $i$.  Furthermore, if $j'$ is an integer such that $jj' \equiv 1 \modu r)$, then $\alpha_{j,k,z} (s^{j'})= s$.  Finally, if $k'$ is an integer such that $(j+nk)k' \equiv -j'k \modu r/p)$, then 
\[ \alpha_{j,k,z} \(c^{pk'} t^{pj'}\) = c^{p(j+nk)k'}\cdot c^{pj'k} t^{pjj'} = t^p.\]  Since there exist such integers $j',k'$ by assumption and  since $s_1,\dots s_{n-1},s,t^p$ generate $G(r,p,n)$, it follows that our map is surjective and hence an automorphism.
\end{proof}

\def\id{\mathrm{id}}

Given $g \in G(r,1,n)$, let $\Ad(g) : x \mapsto gxg^{-1}$ denote the corresponding inner automorphism.  Each such $\Ad(g)$ of course restricts to an automorphism of the normal subgroup $G(r,p,n)$, and with slight abuse of notation we regard $\Ad(g)$ for $g\in G(r,1,n)$ as an element of $\Aut\(G(r,p,n)\)$.  The following lemma gives a useful characterization of which maps $\Ad(g)$ restrict to elements of $\Inn\(G(r,p,n)\)$.

\begin{lemma}\label{inner}
Let $r,p,n$ be positive integers with $p$ dividing $r$.  If $g \in G(r,1,n)$, then the following are equivalent:
\begin{enumerate}
\item[(i)] $\Ad(g)$ restricts to an inner automorphism of $G(r,p,n)$.
\item[(ii)] $\Ad(g)(\pi)$ is conjugate to $\pi$ in $G(r,p,n)$ for all $\pi \in S_n$.
\item[(iii)] $\Delta(g) \in d \ZZ_r$ where $d = \gcd(p,n)$.

\end{enumerate}
\end{lemma}

\begin{proof}
The lemma is trivially true if $n=1$ so assume $n\geq 2$.  
Clearly (i) implies (ii), so assume (ii) holds.  Choose $a \in [0,p-1]$ such that  $g = g' t^a$ for some $g' \in G(r,p,n)$ and let $\pi = (1\cyc 2\cyc \cdots\ n)^{-1} \in S_n$.   Then (ii) implies 
\[\Ad(t^a)(\pi)  =  \(-ae_1 + ae_2,\pi \)= (x,\sigma) \pi (x,\sigma)^{-1}\] for some $(x,\sigma) \in G(r,p,n)$.  Conjugating both sides of this equation by $\sigma^{-1}$ gives 
\be\label{ref} \(-a e_{\sigma^{-1}(1)} + a e_{\sigma^{-1}(2)}, \sigma^{-1} \pi \sigma\) = \( (x_{n}-x_1) e_1 + (x_{1}-x_2) e_2 + \dots + (x_{n-1} - x_n) e_n,\pi\).\ee  Since $C_{S_n}(\pi) = \langle \pi \rangle$, we must have $\sigma = \pi^{i-2}$ for some $i \in [1,n]$.  In this case  $\sigma^{-1}(2)=i$ and $\sigma^{-1}(1) = i-1$ or $n$, and so if we make the abusive definition $x_{j-n} = x_j$ for $j\in [1,n]$, then equation (\ref{ref}) implies 
 \[ x_{i-1} - x_{i} = a,\qquad  x_{i-2}-x_{i-1}=-a,\qquad\text{and}\qquad x_{i-2} = x_{i-3} = \dots = x_{i-n} \overset{\mathrm{def}}= b \in \ZZ_r.\]  From these identities, one computes $a + bn = x_1 + x_2 + \dots + x_n \in p \ZZ_r$, so we have $\Delta(g) +bn = \Delta(g') + a +bn \in p \ZZ_r \subset d \ZZ_r$.  Since $bn \in d\ZZ_r$, (iii) follows.  

Finally assume (iii) holds.  Then since $d = ip + jn$ for some $i,j \in \ZZ$, it follows that (viewing $\ZZ_r$ as a $\ZZ$-module) $\Delta(g) - kn \in p\ZZ_r$ for some $k\in \ZZ$ in which case $gc^{-k} \in G(r,p,n)$.  Since $c$ is central, we have $\Ad(g) = \Ad(gc^{-k})$ which gives (i).  
 %
 %
 %
 %
 %
%
\end{proof}

In almost all cases every automorphism of $G(r,p,n)$ arises by  composing $\Ad(g)$ for some $g \in G(r,1,n)$ with some $\alpha_{j,k,z}$.  To make this precise, we first require the next lemma.

\begin{lemma}\label{main-aut}
Let $r,p,n$ be positive integers with $p$ dividing $r$. Then every automorphism of $G(r,p,n)$ which preserves the normal subgroup $(\ZZ_r)^n \cap G(r,p,n) = \{ g \in G(r,p,n) : |g| = 1\}$ is of the form 
\[\Ad(g) \circ \alpha_{j,k,z},\qquad\text{for some $g \in G(r,1,n)$ and $j,k,z$ as in (\ref{aut-prop-cond}),}\] unless $(r,p,n)$ is $(1,1,6)$.
%

    \end{lemma}

The following elementary proof uses many of the same arguments as the proof of Proposition 4.1 in \cite{MM}, which describes a related but less specific result.

\def\proj{\mathrm{proj}}
\begin{proof}
Let $G=G(r,p,n)$ and $N = (\ZZ_r)^n \cap G$, so that $G = N \rtimes S_n$. 
For the present we assume $n\neq 6$.  Fix $\upsilon \in \Aut(G)$ and suppose $\upsilon(N) = N$.    
If $\proj : G \to S_n$ denotes the homomorphism $(x,\pi ) \mapsto \pi$, then it follows that  map $\pi \mapsto \proj \circ \upsilon(0,\pi)$ is an inner automorphism of $S_n$, and so  $\upsilon = \Ad(\omega) \circ \upsilon'$ for some $\omega \in S_n$ and some $\upsilon' \in \Aut(G)$ with 
 \be \label{form} \upsilon'(s_i) = \(  x_{i,1} e_1 + \cdots + x_{i,n} e_n ,(i\cyc i+1)\)\,\qquad\text{for some choice of }x_{i,j} \in \ZZ_r\ee for $i=1,\dots,n-1$.  Since $s_i^2=1$ we must have $x_{i,i} = -x_{i,i+1}$ and $2x_{i,j} = 0$ for $j\notin \{i,i+1\}$.  
 Since $s_1$ and $s_j$ commute for $j>2$, inspection of the equal expressions $\upsilon'(s_1)$ and $\upsilon'(s_j) \upsilon'(s_1) \upsilon'(s_j)$ shows that $x_{1,3} =x_{1,4}=\dots =  x_{1,n}$.  Therefore  
 \[ \upsilon'(s_1) = z s^{a_1} s_1,\qquad\text{for some $a_1 \in \ZZ_r$ and some $z \in C(r,p,n)$ with $z^2 = 1$.}\]  The conjugacy class of $z s^{a_1} s_1$ consists of elements of the form 
$z \( a e_i - a e_j, (i\cyc j)\)$ for $a \in \ZZ_r$ and $1\leq i<j\leq n$.  The element $\upsilon'(s_i)$ must be of this form, as well of the form (\ref{form}), so we conclude that 
$\upsilon'(s_i) = z \( a_i e_i - a_i e_{i+1}, (i\cyc i+1)\)$ for some $a_i \in \ZZ_r$ for each $i=1,\dots,n-1$.  
 %
 Once can check that if 
\[ \barr{c} y =\(  \sum_{i=1}^n \sum_{j=i}^n a_je_i, 1\) \in G(r,1,n)\qquad\text{then}\qquad y^{-1} \cdot \upsilon'(s_i)\cdot  y = z s_i\text{ for all $i$},\earr\] and so $\upsilon = \Ad({\omega y}) \circ \upsilon''$ where $\upsilon'' \in \Aut(G)$ has $\upsilon''(s_i) = z s_i$ for all $i$.  Since $N$ is normal in $G(r,1,n)$, it follows that $\upsilon''(N) = N$.

Since 
$s$ and $t$ commute with $s_j$ for $j>2$ and $j>1$ respectively, 
and since $s_1 s s_1 = s^{-1}$, 
it follows that we can write 
\[ \upsilon''(s) = z' s^{j}\quad\text{and}\quad\upsilon''(t) = z'' t^{pj'},\qquad\text{for some $j,j' \in \ZZ_r$ and  $z' ,z''\in C(r,p,n)$ with $(z')^2=1$.}\]
If $n=2$ then it follows that either  $z'=1$ or $z' = s^{r/2}$, and so in this case we lose no generality by assuming $z'=1$.  If $n>2$, then since $s_1 s^j = (s_2 s^j s_2)^{-1} s_1 ( s_2 s^j s_2)$, we have 
 \[ zz' s_1 s^j = \upsilon''(s_1s) = \upsilon''((s_2 s s_2)^{-1} s_1 ( s_2 s s_2)) = z s_1 s^j\] so $z'= 1$ automatically.  Hence $\upsilon''(s) = s^j$ for some $j \in \ZZ_r$.  
Since $s^p = t^p s_1 t^{-p} s_1$, we obtain 
\[  \( pj e_1 - pj e_2,1\) = \upsilon''(s^p) = \upsilon''(t^p)\cdot  zs_1 \cdot \upsilon''(t^p)^{-1}\cdot  zs_1 = \( pj' e_1 -pj'e_2, 1\).\]  Therefore $t^{pj} =t^{pj'}$ so we can assume $j'=j$.  Since $t^p$ has order $r/p$, the central element $z''$ must be of the form $c^{pk}$ for some integer $k$, and so $\upsilon''(t^p) = c^{pk} t^{pj}$.  But now $\upsilon''$ agrees with the map $\alpha_{j,k,z}$ on the generators $s_1,\dots,s_{n-1},s,t^p$, and so we conclude by Proposition \ref{aut-prop} that $\upsilon'' = \alpha_{j,k,z}$.  Thus $\upsilon = \Ad(\omega y) \circ \alpha_{j,k,z}$ as desired.

To finish our proof we must treat the case $n=6$ and $r>1$.  In this situation, $N$ is a characteristic subgroup by \cite[Lemma 4.2]{MM} and so the map $\pi \mapsto \proj\circ \upsilon( 0,\pi)$ again induces an automorphism of $S_n$.  Our desired conclusion will follow if we can show that this automorphism is inner, since then we can invoke all of the preceding arguments.    This is shown in the last paragraph of the proof of Proposition 4.1 in \cite{MM}.
\end{proof}

Barring a finite number of cases, the subgroup $(\ZZ_r)^n \cap G(r,p,n)$ is typically characteristic and so every automorphism of $G(r,p,n)$ is of the form given in the lemma. To account for the possible exceptions, we define additional automorphisms $\eta_{r,p,n}, \eta_{r,p,n}' \in \Aut\(G(r,p,n)\)$ on generators by
\[ \ba \eta_{2,1,2} :  s_1 & \mapsto  t \\ 
t & \mapsto s_1\ea
\qquad\qquad
 \ba \eta_{2,2,2} :  s_1 & \mapsto s  \\ s_1' & \mapsto s_1\\ s & \mapsto  s_1' \ea
\qquad\qquad
   \ba \eta_{4,2,2}: s_1 & \mapsto  t^2\\ s_1' & \mapsto  s_1 \\ t ^2& \mapsto s_1' \ea 
\qquad\qquad
\ba \eta_{3,3,3} :s_1 & \mapsto s_2 \\ s_2 & \mapsto  s_1' \\ s_1' & \mapsto s_1
\ea 
\]

\[ \ba \eta_{3,3,3}' :s_1 & \mapsto s_1 \\ s_2 & \mapsto  s_2 \\ s_1' & \mapsto s_2'
\ea 
\qquad\qquad
  \ba \eta_{2,2,4}: s_1 & \mapsto  s_1' \\ s_2 & \mapsto  s_2 \\ s_3 & \mapsto s_1 \\ s_1' & \mapsto s_3 \ea 
\qquad\qquad
 \ba \eta_{1,1,6} : 
 s_1 & \mapsto  (1\cyc 2)(3\cyc 4)(5\cyc 6) \\ 
 s_2 & \mapsto  (1\cyc 5)(2\cyc 3)(4\cyc 6) \\ 
 s_3 & \mapsto (1\cyc 2)(3\cyc 6)(4\cyc 5) \\
 s_4 & \mapsto  (1\cyc 5)(2\cyc 6)(3\cyc 4) \\ 
 s_5 & \mapsto  (1\cyc 2)(3\cyc 5)(4\cyc 6)  \ea
 \]
 and in all other cases $\eta_{r,p,n} = 1$ and $\eta_{r,p,n}' = 1$.  Thus $\eta_{r,p,n}'$ is the identity unless $(r,p,n)=(3,3,3)$.  
 
 Many of these automorphisms are well known: for example, $\eta_{2,1,2}$ and $\eta_{2,2,4}$ are the graph automorphisms $^2B_2$ and $^3D_4$, and $\eta_{1,1,6}$ is the outer automorphism of $S_6$.  The automorphisms $\eta_{4,2,2}$ and $\eta_{3,3,3}$, $\eta_{3,3,3}'$ and $\eta_{2,2,4}$ derive from normal embeddings of $G(4,2,2)\vartriangleleft G_6$ and $G(3,3,3) \vartriangleleft G_{26}$ and $G(2,2,4)\vartriangleleft G_{28}$ in exceptional groups.  For more information on the structure of this embedding, see \cite[Proposition 3.13]{BMM} or \cite[Section 3]{MM}.  



\begin{theorem}
Let $r,p,n$ be positive integers with $p$ dividing $r$.  
Then every automorphism of $G(r,p,n)$ is of the form 
\[ (\eta_{r,p,n})^{i_1}\circ (\eta_{r,p,n}')^{i_2}\circ \Ad(g) \circ \alpha_{j,k,z} \] for some $g \in G(r,1,n)$, $i_1,i_2,j,k \in \ZZ$, and $z \in C(r,p,n)$.
\end{theorem}

\begin{remark}
Of course, by Lemma \ref{aut-prop} the given expression is an automorphism of $G(r,p,n)$ if and only if $\gcd(j,r) = \gcd(j+nk,r/p)=1$ and $z^2=1$.  Note, furthermore, that $\eta_{r,p,n} =\eta'_{r,p,n} = 1$ if $(r,p,n)$ is not $(2,1,2)$, $(2,2,2)$, $(4,2,2)$, $(3,3,3)$, $(2,2,4)$, or $(1,1,6)$.  In these cases, it is easy to see that the reflection-preserving automorphisms of $G(r,p,n)$ are precisely the maps of the form $\Ad(g)\circ \alpha_{j,0,1}$, as predicted by \cite[Theorem 7.1]{ShiWang}.
\end{remark}

\begin{proof}
One can check directly that the theorem holds if $(r,p,n)$ is $(2,1,2)$, $(2,2,2)$, $(4,2,2)$, $(3,3,3)$, $(2,2,4)$, or $(1,1,6)$; we have done so using the computer algebra system {\tt{GAP}}.  If $(r,p,n)$ is not one of these exceptions, then by  Lemma 4.2 in \cite{MM} the subgroup $(\ZZ_r)^n \cap G(r,p,n)$ is  characteristic, and hence preserved by $\Aut(G)$, in which case the theorem follows immediately from Lemma \ref{main-aut}.   
\end{proof}

We can be a little more specific about the uniqueness of the decomposition given in the theorem, and this will allow us to give a formula for the order of $\Aut\(G(r,p,n)\)$.  
Given integers $j,k$ and $z \in C(r,1,n)$, we adopt the shorthand 
\[ \beta_j \overset{\mathrm{def}} =\alpha_{j,0,1} : (x,\pi)\mapsto (jx,\pi)\qquad\text{and}\qquad \gamma_{k,z} \overset{\mathrm{def}}= \alpha_{1,k,z}: (x,\pi) \mapsto z^{\ell(\pi)} c^{\Delta(x)\cdot k} (x,\pi).\] One checks that 
\[ \beta_j \circ\beta_{j'} 
= 
\beta_{jj'}, \qquad
 \gamma_{k,z} \circ \gamma_{k',z'} 
 = \gamma_{k'',zz'} , \qquad\text{and}\qquad
\beta_j \circ \gamma_{k,z} = \gamma_{k,z} \circ \beta_j 
 = \alpha_{j,jk,z},\] 
where $k'' = k+k' + nkk'$.  Since in this notation $\alpha_{0,1,1} =  \beta_1 = \gamma_{0,1} $ all equal  the identity automorphism of $G(r,p,n)$, it follows that the sets  
\[  B \overset{\mathrm{def}}= \{ \beta_j : \gcd(j,r)=1\}  \qquad\text{and}\qquad C\overset{\mathrm{def}} = \{\gamma_{k,z}  : \gcd(1+nk,r/p)=1,\ z \in C(r,p,n),\ z^2=1\} 
\] 
are subgroups of $ \Aut\(G(r,p,n)\)$.  Define additionally the subgroups 
\[
X\overset{\mathrm{def}} = \langle \eta_{r,p,n}, \eta_{r,p,n}' \rangle\qquad\text{and}\qquad A \overset{\mathrm{def}}= \{ \Ad(g) : g \in G(r,1,n) \}\] of $\Aut\(G(r,p,n)\)$.  
Using the previous result, we have the following.  

\begin{proposition}\label{before}
Let $r,p,n$ be positive integers with $p$ dividing $r$.  Then 
\begin{enumerate}
\item[(i)]  $\Aut\(G(r,p,n)\) = XABC$ and $XA $ is a normal subgroup of $\Aut\(G(r,p,n)\)$.

\item[(ii)] \begin{enumerate}
\item[a.] If $n=1$ then $B=C$, and if $n>1$ then $B\cap C = \{1\}$. 
\item[b.] If $n=2$ then $A\cap BC$ is the set of all automorphisms of the form $\alpha_{j,k,z}$ with $(j,k) = (-1,1)$ or $(1,0)$, 
and if $n\neq 2$ then $A \cap BC = \{1\}$. 
\item[c.] If $(r,p,n) = (3,3,3)$ then $X\cap ABC = A$ and if $(r,p,n)\neq (3,3,3)$ then 
$X \cap ABC = \{1\}$.
 \end{enumerate}
 
\item[(iii)] If $n>2$ then $\Aut\(G(r,p,n)\) \cong (XA\rtimes B) \times C$.
\end{enumerate}
\end{proposition}

This result is closely related to Proposition 4.1 in \cite{MM}, which asserts that any automorphism of $G(r,p,n)$ is the composition of an automorphism which preserves the pseudo-reflections in $G(r,p,n)$ and a central automorphism.  The subgroup $C\subset \Aut\(G(r,p,n)\)$ is the set of central automorphisms of $G(r,p,n)$ and $XAB\subset \Aut\(G(r,p,n)\)$ is the subgroup of automorphisms which preserve the reflections.  Furthermore, if we view $G(r,p,n)$ as a subgroup of $\GL(\CC^n)$ whose elements are generalized permutation matrices, then the subgroup $XA$ consists of all automorphisms induced from elements of the normalizer of $G(r,p,n)$ in $\GL(\CC^n)$.  The subgroup $B$ likewise consists of all automorphisms induced from elements of the Galois group of $\QQ(\zeta_r)$ over $\QQ$.

\begin{proof} 
To prove the first half of (i), one checks that $BC$ is precisely the set of all automorphisms of $G(r,p,n)$ of the form $\alpha_{j,k,z}$; hence, by the preceding theorem $\Aut\(G(r,p,n)\) = XABC$.  One can check the remaining assertions in the finite number of cases when $X \neq \{1\}$ directly, by hand or with a computer algebra system.  Therefore assume $X =\{1\}$. 

 
   To show $A \vartriangleleft \Aut\(G(r,p,n)\)$ we observe  that if $g = (y,\omega) \in G(r,1,n)$ then $\alpha_{j,k,z} \circ \Ad(g) \circ \alpha_{j,k,z}^{-1} = \Ad({g'})$ where $g' = (jy,\omega)$.  
Our description of $B\cap C$ is trivially verified.  Suppose $n=2$ and $z \in C(r,p,n)$ has $z^2=1$.  We have two cases: either $z=1$ or $z = \( \frac{r}{2} e_1 + \frac{r}{2} e_2,1\)$, the latter of which can only  occur if $r$ is even.  One checks that $\alpha_{-1,1,z} = \Ad({g})$ for $g = s_1$ if $z=1$ and for $g=\( \frac{r}{2} e_1, (1\cyc 2)\)$ if $z\neq 1$.  Likewise, $\gamma_{0,z}$ is the identity if $z=1$ and $\Ad(g)$ for $g = \( \frac{r}{2} e_1,1\)$ if $z\neq 1$.  Since $n=2$ we have \[ \Ad(g)(s_1 ) = s^a s_1\text{ for some $a \in \ZZ$}\qquad\text{and}\qquad \ba \Ad(g)(s) &= s\\ \Ad(g)(t^p) &= t^p\ea \quad \text{ or }\quad \ba \Ad(g)(s) &= s^{-1} \\ 
 \Ad(g)(t^p) &= c^pt^{-p}\ea\] and it follows that only elements of the form $\alpha_{j,k,z}$ with $(j,k)=(1,0)$ or $(-1,1)$ can be contained in $A \cap BC$, which proves (ii).

If $n=1$ then $A = \{1\}$.  Suppose $n>2$ and $\Ad(g) = \alpha_{j,k,z}$ for some $g = (y,\pi) \in G(r,1,n)$ and $j,k,z$.  This implies $\pi \sigma \pi^{-1} = \sigma$ for all $\sigma \in S_n$ so $\pi = 1$ since $Z(S_n)$ is trivial.  Consequently $\Ad(g)$ fixes both $s$ and $t$ so we must have $j=1$ and $k=0$. Furthermore, $\Ad(g)(s_1) = s^{y_1-y_2} s_1$ so we must have $y_1=y_2$ and $z=1$ since $n>2$.  Therefore $\alpha_{j,k,z} = 1$ so $A\cap BC = \{1\}$.  Since  each element in $C$ commutes with all elements of $A$ and $B$, part (iii) now follows from (i) and (ii).
\end{proof}

From this proposition we are now able to derive a formula for the order of the automorphism group of $G(r,p,n)$.  Here $\phi(x)$ denotes Euler's totient function, which we recall is defined as the number of positive integers $y \leq x$ with $\gcd(x,y)=1$.  

\begin{corollary}
Let $r,p,n$ be positive integers with $p$ dividing $r$.  Assume $n>1$ and write $e$ for the greatest divisor of $r/p$ with $\gcd(e,n) =1$.  Then 
\begin{enumerate}
\item[] $| \Aut\(G(r,p,n)\)| = \displaystyle \frac{c_{r,p,n}}{c'_{r,p,n}} \cdot \phi(r) \cdot \phi(e)/e \cdot n!\cdot  r^n/p$,
\end{enumerate}
and
\begin{enumerate}
\item[] $|\Out\(G(r,p,n)\)| = c_{r,p,n}   \cdot \phi(r)\cdot  \phi(e )/e \cdot r/p \cdot \gcd(p,n)$,

\item[]$ |Z\(G(r,p,n)\)| = c'_{r,p,n}\cdot  r/p \cdot \gcd(p,n)$,

\end{enumerate}
 where 
 \[ \ba &c_{1,1,2} = 1,\quad &&c_{2,2,2} =3,\quad &&c_{2,1,2} =1,\quad &&c_{4,2,2} = 3/2,\quad &&c_{3,3,3} = 4,\quad &&c_{2,2,4} = 6,\ \quad &&c_{1,1,6} = 2,\\
 &c'_{1,1,2} =2, \quad&& c'_{2,2,2} =2,\ea\] and in all other cases
\[ c_{r,p,n} = \left\{\ba
&1/2, &&\quad\text{if $n=2$}, \\
&1, && \quad\text{if $r$ is odd and $n>2$}, \\
&1, & &\quad\text{if $p$ is even but $r/p$ and $n>2$ are odd}, \\ 
&2,&&\quad\text{otherwise},
\ea\right.
\qquad \text{and}\qquad
c'_{r,p,n} = 1.
\]

\end{corollary}

\begin{remark} If $n=1$ then $G(r,p,n) \cong \ZZ_{r/p}$ so $|\Aut\(G(r,p,n)\)| = |\Out\(G(r,p,n)\)| = \phi(r/p)$.
\end{remark}

\begin{proof}
The formula for the order of the center of $G=G(r,p,n)$ follows immediately from Lemma \ref{center}, so it suffices to prove our formula for $|\Aut(G)|$.   If $(r,p,n)$ is one of the exceptional cases $(1,1,2)$, $(2,2,2)$, $(2,1,2)$, $(4,2,2)$, $(3,3,3)$, $(2,2,4)$, or $(1,1,6)$, then our formula asserts that $|\Aut(G)|$ is 1, 6, 8, 48, 432, 1152, or 1440, respectively.  One easily checks that these orders are correct: all but two of the exceptions are Weyl groups whose automorphisms are well known (e.g., see \cite{B1969}), and one can compute the outer automorphisms of $G(4,2,2)$ and $G(3,3,3)$ by hand or with a computer.

Assuming $(r,p,n)$ is not one of these exceptions, we have  $|A| = |G(r,1,n)| / |C(r,1,n)| = n!\cdot r^{n-1}$  by Lemma \ref{center} and $|B| = \phi(r)$.  To compute $|C|$, we note that $\gamma_{k,z} = \gamma_{k',z'}$ if and only if $k\equiv k' \modu r/p)$ and $z=z'$.  Thus the elements of $C$ are in bijection with all choices of $k \in [0,r/p-1]$ and $z \in C(r,p,n)$ such that $\gcd(1+nk,r/p) =1$ and $z^2=1$.  To satisfy these conditions, the central element $z$ must be the identity if $r$ is odd, or if  $r$ is even but $nr/2$ is not a multiple of $p$, which occurs if and only if $p$ is even but $r/p$ and $n$ are odd.  In all other cases, $z$ can be either $1$ or $c^{r/2}$.  Additionally, since $1+nk$ is coprime to $r/p$ if and only if $1+nk$ is coprime to $e$,  it follows that there are $\phi(e) \cdot (r/p)/e$ possible choices of $k \in [0,r/p-1]$ with $\gcd(1+nk,r/p)=1$.  Thus 
$|C| = \wt c_{r,p,n} \cdot \phi(e)/e \cdot r/p$ where 
\[\wt c_{r,p,n}  = \left\{\ba
&1, && \quad\text{if $r$ is odd, or if $p$ is even but $r/p$ and $n$ are odd}, \\ 
&2,&&\quad\text{otherwise}.
\ea\right.\]  As $n>1$ and $X = \{1\}$, by the preceding proposition we thus have 
\[|\Aut(G)| =\frac{ |A| |B| |C|}{|A\cap BC| |B\cap C|}=\frac{\wt c_{r,p,n}}{|A\cap BC|} \cdot \phi(r) \cdot \phi(e)/e \cdot n!\cdot  r^n/p.\]   If $n>2$, then $|A\cap BC| = 1$ and if $n=2$ then it follows from part (ii) of Proposition \ref{before} that $|A\cap BC| = 2\wt c_{r,p,n} $, so in both cases we obtain $\wt c_{r,p,n}/|A\cap BC| = c_{r,p,n} / c'_{r,p,n}$ as desired. 
\end{proof}


\section{Classifying Generalized Involution Models}\label{classifying-sec}

In this final section we complement the results of Section \ref{positive} by proving that $G(r,p,n)$ does not have a generalized involution model if $\gcd(p,n) \geq2$, unless $\gcd(p,n)=n=2$ and $r/p$ is odd.
Combining this with a computer assisted investigation of the exceptional groups, we will be able to completely determine which finite complex reflection groups have generalized involution models.

\subsection{Observations and Reductions}

We begin by observing that finding all generalized involution models of a group often is equivalent to classifying the generalized involution models defined with respect to a single, fixed automorphism.
Say that an automorphism $\tau \in \Aut(G)$ of a group $G$ is \emph{class-preserving} if $\leftexp\tau g$ is conjugate to $g$ for all $g \in G$, or equivalently if $\psi \circ \tau = \psi$ for all $\psi \in \Irr(G)$.  
Clearly all inner automorphisms are class-preserving, but a finite group can possess outer automorphisms which are class-preserving, as was first shown by  Burnside \cite{Burnside}.
The non-existence of class-preserving outer automorphisms can greatly reduce the problem of finding all generalized involution models of a group by the following lemma.

\begin{lemma}\label{classifying}
Let $G$ be a finite group with an automorphism $\tau \in \Aut(G)$  such that $\tau^2=1$ and \be\label{i}\sum_{\psi \in \Irr(G)} \psi(1) = | \{ g \in G : g\cdot \leftexp{\tau}g = 1\}|.\ee  If $G$ has no class-preserving outer automorphisms, then the following hold:
\begin{enumerate}
\item[(i)] The image of $\tau$ in $\Out(G)$ is central.

\item[(ii)] Any generalized involution model for $G$ can be defined with respect to $\tau$.    
\end{enumerate}
\end{lemma}

\begin{remark}
The conclusion of the lemma can fail if $G$ has class-preserving outer automorphisms.  Wall \cite{Wall}  showed that the semidirect product  $G=  \ZZ_m \rtimes \ZZ_m^\times$ consisting of all pairs $(a,x) \in \ZZ_m \times \ZZ_m^\times$ with the multiplication 
\[ (a,x)(b,y) = (a + xb, xy),\qquad\text{for }a,b \in \ZZ_m,\ x,y \in \ZZ_m^\times\]
 has a class-preserving outer automorphism $\tau$ of order two if $m$ is divisible by 8.  Taking $m=8$ gives a group $G$ of order 32, the smallest group with a class-preserving outer automorphism.  One can check using a computer algebra system (we used {\tt{GAP}}) that this $G$ has a generalized involution model with respect to $1 \in \Aut(G)$ but not with respect to the class-preserving outer automorphism $\tau$, even though (\ref{i}) holds.
\end{remark}

\begin{proof}
Assume $G$ has no class-preserving outer automorphisms.  If $\alpha \in \Aut(G)$ and $\tau'=\alpha\circ \tau \circ \alpha^{-1}$, then by Theorem \ref{bg-thm}, we have $\epsilon_\tau(\psi) =1$ and $\epsilon_{\tau'}(\psi) = \epsilon_\tau(\psi\circ \alpha) =1$ for all $\psi \in \Irr(G)$.  Therefore, by Proposition 2 in \cite{BG2004}, $\leftexp{\tau}g$ is conjugate to $\leftexp{\tau'}g$ for all $g \in G$, so $\alpha\circ \tau \circ \alpha^{-1} \circ \tau^{-1}$ is class-preserving, and therefore an inner automorphism.  This proves (i).

Now suppose $G$ has a generalized involution model with respect to $\upsilon \in \Aut(G)$ with $\upsilon^2=1$.  
By Theorem \ref{bg-thm}, it follows that $\epsilon_\tau(\psi) = \epsilon_\upsilon(\psi) = 1$ for all $\psi \in \Irr(G)$, so by \cite[Proposition 2]{BG2004} each $g \in G$ is conjugate to both $\leftexp{\tau}g^{-1}$ and $\leftexp{\upsilon}g^{-1}$.  Replacing $g$ with $g^{-1}$, one sees that $\leftexp{\tau}g$ is therefore conjugate to $\leftexp{\upsilon}g$ for all $g \in G$, which suffices to show that $\tau\circ \upsilon^{-1} $ is a class-preserving automorphism.  Hence $\tau = \Ad(x) \circ \upsilon$  for some $x \in G$.  Since $\tau^2=\upsilon^2=1$, the element $z \overset{\mathrm{def}} = x\cdot \leftexp{\upsilon}x \in G$ is central.  Fix  $\psi \in \Irr(G)$ and let $\omega_\psi(z) \overset{\mathrm{def}}= \frac{\psi(z)}{\psi(1)}$ denote the value of its central character at $z$; then $\psi(zg) = \omega_\psi(z) \psi(g)$ for all $g\in G$, and it follows that 
 $\omega_\psi(z) \epsilon_\tau(\psi) = \epsilon_\upsilon(\psi) = \epsilon_\tau(\psi)= 1$ so $\psi(z) = \psi(1)$.  Since this holds for all irreducible characters of $G$, we have $z=1$.  This means that $g\cdot\leftexp{\tau}g = gx\cdot \leftexp{\upsilon}g\cdot  x^{-1} = (gx)\cdot \leftexp{\upsilon}(gx)$, and it follows that the map $\cI_{G,\tau} \to \cI_{G,\upsilon}$ given by $g \mapsto gx$ is an isomorphism of $G$-sets.  In particular, the twisted conjugacy classes with respect to $\tau$ and $\upsilon$ are in bijection and have the same twisted centralizers.  Therefore the generalized involution model with respect to $\upsilon$ can also be defined with respect to $\tau$, which proves (ii).
\end{proof}

As a consequence of this result, to determine whether a group $G$ with no class-preserving outer automorphisms has a generalized involution model, one only needs to check (1) if there exists $\tau \in \Aut(G)$ with $\tau^2=1$ such that (\ref{i}) holds,
and (2) if $G$ has a generalized involution model with respect to $\tau$.  
This strategy is especially apposite for irreducible complex reflection groups in light of the following result given in a slightly different form as Proposition 3.1 in \cite{MM}.

\def\Out{\mathrm{Out}}

\begin{lemma}\label{cmplx} (Marin, Michel \cite{MM}) A finite complex reflection group has no class-preserving outer automorphisms.
\end{lemma}

\begin{remark}
We remark that it is a tedious but not overly difficult exercise to prove the lemma directly for the irreducible groups $G(r,p,n)$, and via computer calculations for the exceptional groups.  
The lemma then holds for all finite complex reflection groups because a class-preserving automorphism of a direct product must restrict to a class-preserving automorphism of each factor.
\end{remark}


These results become especially useful when doing calculations. To determine which of the exceptional irreducible complex reflection groups $G_4,\dots,G_{37}$ have generalized involution models, we will rely on a computer-assisted brute force search.  The preceding lemmas greatly diminish the size of this calculation, because they show that one needs to examine at most one automorphism for each group to determine if a generalized involution model exists.  In Table \ref{tbl1} we provide a list of automorphisms $\tau \in \Aut(G_i)$ for which (\ref{i}) holds, if this is possible.  These automorphisms are defined on the generators $s,t,u,v,w$ which appear in the presentations for $G_4,\dots,G_{37}$ in the appendix of \cite{BMR}.  These generators coincide with the generators for the exceptional groups in the {\tt{GAP}} package {\tt{CHEVIE}}, which allows one easily to compute things with this data.

\begin{table}[h]
\[
\barr{|c|c|} \hline & \\[-10pt]   \text{Exceptional Group $G$} & \text{\ $\tau \in \Aut(G)$ with $\tau^2=1$ and $\epsilon_\tau(\psi)=1$ for all $\psi \in \Irr(G)$}
  \\[-10pt] 
& 
  \\\hline
& 
\\[-10pt]  
    \barr{c}       G_4,\
           G_5,\
           G_6,\
           G_8,\
           G_9,\
           G_{10}, \\
           G_{14},\
           G_{16},\
           G_{17},\
           G_{18},\
           G_{20},\
           G_{21} \earr & (s,t) \mapsto \(s^{-1}, t^{-1}\)
  \\[-10pt]
& \\\hline 
& \\[-10pt]
    \barr{c}       G_7,\ G_{11},\ 
           G_{19} \earr
         & (s,t,u) \mapsto \(s,t^{-1},su^{-1}s\)
                  \\[-10pt]
& \\\hline 
& \\[-10pt]
    \barr{c}       G_{12},\ G_{22},\ G_{24},\ G_{25} \earr & (s,t,u) \mapsto \(u^{-1}, t^{-1}, s^{-1}\)
                  \\[-10pt]
& \\\hline 
& \\[-10pt]
    \barr{c}       G_{13} \earr & (s,t,u) \mapsto \(s,u,t\)
                  \\[-10pt]
& \\\hline 
& \\[-10pt]
    \barr{c}       G_{15} \earr & (s,t,u) \mapsto \(s,t,tu^{-1}t\)
 \\[-10pt]&
            \\\hline 
& \\[-10pt]
    \barr{c}       G_{23},\ G_{28},\
           G_{30},\
           G_{35},\ G_{36},\ G_{37} \earr & \text{Identity automorphism}
\\[-10pt]&
    \\\hline 
& \\[-10pt]
\barr{c}    G_{26} \earr & (s,t,u) \mapsto \(s^{-1}, t^{-1}, u^{-1}\)
                              \\[-10pt]&
            \\\hline 
& \\[-10pt]
    \barr{c}       G_{27},\
           G_{29},\ G_{34} \earr & \text{No such $\tau$ exists}
                                   \\[-10pt]&
            \\\hline 
& \\[-10pt]
    \barr{c}      G_{31} \earr & (s,t,u,v,w) \mapsto\( u,t,s,w,v\)
                                   \\[-10pt]&
            \\\hline 
& \\[-10pt]
    \barr{c}       G_{32} \earr & (s,t,u,v) \mapsto \(v^{-1}, u^{-1}, t^{-1}, s^{-1}\)
                                       \\[-10pt]&
            \\\hline 
& \\[-10pt]
    \barr{c}       G_{33} \earr & (s,t,u,v,w) \mapsto \(v,u,t,s,w\)
  \\[-10pt] 
& 
  \\\hline
\earr\]
\caption{Automorphisms $\tau \in \Aut(G_i)$ satisfying the conditions of Theorem \ref{bg-thm}, defined in terms of the presentations given in \cite{BMR}}
\label{tbl1}
\end{table}

Vinroot, elaborating upon the work of Baddeley \cite{B91}, describes in \cite{V} all finite Coxeter groups with involution models in the classical sense; in particular, the only irreducible finite Coxeter groups which fail to have involution models are those of type  $D_{2n}$ ($n>1$), $E_6$, $E_7$, $E_8$, $F_4$, and $H_4$.  If $G$ is a finite Coxeter group then all of its representations are equivalent to real representations, and so by the Frobenius-Schur involution counting theorem, (\ref{i}) holds with $\tau=1$.  Hence, by Lemmas \ref{classifying} and \ref{cmplx}, a finite Coxeter group has a generalized involution model if and only if it has an involution model, and we are left with the following corollary of Theorem 1 in \cite{V}.

\begin{corollary}\label{coxeter-cor}
A finite Coxeter group has a generalized involution model if and only if it has an involution model,  which occurs if and only if
 all of its irreducible factors are of type $A_n$, $B_n$, $D_{2n+1}$, $H_3$, or $I_2(n)$.  
\end{corollary}

\begin{remark} The Coxeter group of type $G_2$ is omitted from this list only because it is isomorphic to the one of type $I_2(6)$.  We note that the Coxeter group of type $I_2(n)$ is the involutory complex reflection group $G(n,n,2)$, and that restricted to this group the map $\tau : g \mapsto \overline g$ is a nontrivial inner automorphism.  Thus, while the group has an involution model in the classical sense, it also has a generalized involution model with respect to $\tau$, which is consistent with  Lemmas \ref{classifying} and \ref{cmplx}.  The same is true for groups of types $A_n$, $B_n$, and $D_{2n+1}$, but vacuously since in these cases the inverse transpose $\tau$ acts as the identity map.
%
\end{remark}

In order to reduce our investigation of finite complex reflection groups to irreducible groups, we require one additional lemma.   This next statement generalizes Lemma 1 in \cite{V} which considers only involution models.

\begin{lemma}\label{products}
If $H_1,\dots, H_n$ are finite groups then $H = \prod_{i=1}^n H_i$ has a generalized involution model if and only if each $H_i$ has a generalized involution model.
\end{lemma}

\begin{proof}
If $H$ has a generalized involution model with respect to $\tau \in \Aut(H)$, then each $h \in H$ is conjugate to $\leftexp{\tau}h^{-1}$  by Theorem \ref{bg-thm} and Proposition 2 in \cite{BG2004}, and so $\tau$ restricts to an automorphism of each factor $H_i$.  Given this fact, it follows that any generalized involution model for $H$ decomposes in an obvious way as a ``product'' of generalized involution models of the factor groups $H_i$, and the proof of the lemma becomes a simple exercise.
 \end{proof}

\subsection{Addressing $G(r,p,n)$ with $\gcd(p,n)=2$}

We now demonstrate that $G(r,p,n)$ does not have a generalized involution model if $\gcd(p,n)=2$, unless $n=2$ and $r/p$ is odd.   Our proof of this proceeds in two steps, and will use somewhat different methods.  We begin in the case when $\gcd(p,n) =2$ and $r/p$ is even.

\begin{lemma}\label{last1}
Let $r,p,n$ be positive integers with $p$ dividing $r$.  If $\gcd(p,n) = 2$ and $r/p$ is even, then $G(r,p,n)$ does not have a  generalized involution model.
\end{lemma}

\begin{proof}
We can tackle this case by much more direct methods than when $r/p$ is odd.  Let $G = G(r,p,n)$ and define $\tau$ as the usual inverse transpose automorphism $g\mapsto \overline g$.  Since $\gcd(p,n)=2$, it follows from Theorem \ref{thm4.5} that equation (\ref{i}) is satisfied, so by Lemmas \ref{classifying} and \ref{cmplx}, we need only show that $G$ has no generalized involution models with respect to $\tau$.
Towards this goal, our strategy is simple.  Since $r/2$ is a multiple of $p$, the central element $z\overset{\mathrm{def}}= c^{r/2} \in G(r,1,n)$ is contained in $G$; here $c$ is defined as in (\ref{gen}).  
We claim that $z$ lies in the commutator subgroup of the twisted centralizer $C_{G,\tau}(\omega)$ for every generalized involution $\omega \in \cI_{G,\tau}$.

  If this holds, then $z$ lies in the kernel of every linear character $\lambda$ of $C_{G,\tau}(\omega)$ and therefore also in the kernel of the induced character $\Ind_{C_{G,\tau}(\omega)}^G(\lambda)$ since $z$ is central.  In this case, if $G$ has a generalized involution model $\{\lambda_i : H_i \to \CC\}$ with respect to $\tau$, then $z$ lies in the kernel of $\sum_i \Ind_{H_i}^G(\lambda_i)$, implying the contradiction  
\[1\neq z \in \bigcap_{\psi \in \Irr(G)} \ker (\psi) = \{1\}.\]

To prove our claim,  suppose $\omega = (a,\pi) \in G$ has $\omega\cdot \leftexp{\tau}\omega = (\pi^{-1}(a)-a,\omega^2) = (0,1) = 1$.  Then $\pi \in S_n$ must be an involution with $\pi(a) =a$.  We lose no generality by conjugating $\omega$ by an element of $S_n \subset G$ since this has the effect of conjugating the twisted centralizer $C_{G,\tau}(\omega)$ and fixing $z$.  Therefore, we can assume that $\omega =  (1\cyc 2)(3\cyc 4)\cdots(2k-1\cyc 2k)$ for some $k\leq n/2$, in which case $\pi(a) = a$ implies  $a_{2i-1} = a_{2i}$ for all $i=1,\dots,k$.  Since $r$ and $p$ are even and $\Delta(a) \in p\ZZ_r$,  the number of $a_i \notin 2\ZZ_r$ is even; therefore, letting $\ell = n/2-k$, there are distinct indices $\{ i_1, j_1,\dots, i_\ell, j_\ell \} = [2k+1,n]$ such that $a_{i_t} - a_{j_t} \in 2\ZZ_r$ for all $t=1,\dots,\ell$.  For each $t$, let $b_t \in \ZZ_r$ such that $2b_t = a_{i_t} - a_{j_t}$.    Now define $g = ( x,\sigma) \in G$ by 
\[ \sigma =(1\cyc 2)\cdots(2k-1\cyc 2k)(i_1\cyc j_1)\cdots(i_\ell\cyc j_\ell) \in S_n
\qquad\text{and}\qquad x_i = \left\{\barr{ll} 0, & \text{if }i \in [1,k], \\ b_t, &\text{if } i = j_t, \\ -b_t,&\text{if }i=i_t.\earr\right.\]  One can check that we then have $\sigma\in C_{S_n}(\pi)$, $\pi(x) = x$, and 
$a+2x = \sigma(a)$, and so 
\[ g\cdot \omega\cdot  \leftexp{\tau}g^{-1} = (\sigma^{-1} \pi^{-1}(x) + \sigma^{-1}(a) + \sigma^{-1}(x), \sigma \pi \sigma^{-1}) 
= \( \sigma^{-1}(a + 2x), \pi\) = \omega.\]  Thus $g \in C_{G,\tau}(\omega)$.  Since $r$ is divisible by 4 and $r/2$ is divisible by $p$, we can define $h = (y,1) \in G$  by setting $y \in (\ZZ_r)^n$ to have 
\[ y_i =  \left\{\barr{ll} r/4, & \text{if }i \in \{1,3,\dots,2k-1\}, \\
 -r/4,&\text{if }i \in \{2,4,\dots,2k\}, \\
r/2, &\text{if } i = i_t, \\ 
0,&\text{if }i=j_t.\earr\right.
\] 
Observe that $\pi^{-1}(y) = -y$ since $r/2=-r/2$, so $h\cdot \omega\cdot \leftexp{\tau}h^{-1} = \omega h^{-1}h = \omega$ and $h \in C_{G,\tau}(\omega)$.  Our claim now follows by calculating $ghg^{-1}h^{-1} = (\sigma^{-1}(y) -y,1) = z$, which completes the proof.
\end{proof}

If $\gcd(p,n) =2$ but $r/p$ is odd, then the crucial step in the preceding proof does not hold.  However, in this case the group $G(r,p,n)$ still fails to have a generalized involution model, provided $n>2$.    To show this, we will use two results from Baddeley's thesis \cite{B91}.

First, recall that a model for a group $G$ is a set $ \{ \lambda_i : H_i \to \CC\}$ of linear characters  of subgroups of $G$ such that 
$ \sum_i \Ind_{H_i}^G(\lambda_i) = \sum_{\psi \in \Irr(G)} \psi$.  Following Baddeley, we say that a model $ \{ \lambda_i : H_i \to \CC\}$ is \emph{based} on a set $\mathcal{S}$ of subgroups of $G$ if for each $i$ there exists a subgroup $H_i' \in  \mathcal{S}$ with a linear character $\lambda_i' : H_i' \to \CC$ such that $\Ind_{H_i}^G(\lambda_i) = \Ind_{H_i'}^G(\lambda_i')$.  Thus $ \{ \lambda_i : H_i \to \CC\}$ is based on the set of subgroups $\{H_i\}$, as well as on any set of subgroups which are conjugate in $G$ to the subgroups $H_i$, as well as on the set of all subgroups of $G$.   For each $n\geq 1$, let $\sG(n)$ denote the set of subgroups of $S_n$ of the form \[W_k \times S_i \times S_j,\qquad\text{where $i,j,k$ are nonnegative integers with $i+j+2k = n$}\] and $W_k \subset S_{2k}$ is the centralizer of the permutation $(1\cyc 2)(3\cyc 4)\cdots(2k-1\cyc 2k) \in S_{2k}$.  The centralizer of any involution in $S_n$ is conjugate to a subgroup of the form $W_k \times S_j = W_k \times S_0 \times S_j \in \sG(n)$ for some $j,k$ with $j+2k = n$, so any involution model for $S_n$ is based on $\sG(n)$.  (This is not a vacuous statement;  \cite{IRS91} constructs an involution model for the symmetric group.)  
Baddeley states the following result as Corollary 4.3.16 in \cite{B91}.

\begin{lemma}\label{baddeley1}
(Baddeley \cite{B91}) Suppose $\cM$ is a model for $S_n$ based on $\sG(n)$.  If $\cM$ contains both the trivial character $\One\in \Irr(S_n)$ and the sign character $\sgn  \in \Irr(S_n)$ then $n=2$.
\end{lemma}

To state our second needed result, let $\Phi : G \to G'$ be a surjective group homomorphism.  Suppose $H\subset G$ is a subgroup and $\psi \in \Irr(H)$.  If $\ker(\psi) \supset \ker(\Phi) \cap H$, then there exists a unique irreducible character $\psi' \in \Irr(\Phi(H))$ such that $\psi = \psi' \circ \Phi$, and we define $\cR_\Phi(\psi) \in \{0 \} \cup \Irr(\Phi(H))$ by
\[ \cR_\Phi(\psi) = \left\{\barr{ll}\psi', &\text{if }\ker (\psi) \supset \ker(\Phi) \cap H, \\ 0,&\text{otherwise.}\earr\right.
\]  
The following appears as Theorem 4.2.3 in \cite{B91}.

\begin{theorem}\label{4.2.3} (Baddeley \cite{B91}) Let $\Phi : G \to G'$ be a surjective  group homomorphism.  If $\cM$ is a model for $G$ and \[\hat \cM = \{ \lambda \in \cM : \cR_\Phi(\lambda) \neq 0\},\] then $\cM_\Phi \overset{\mathrm{def}}=\left\{ \cR_\Phi(\lambda) : \lambda \in \hat \cM \right\}$ is a model for $G'$.
\end{theorem}

We now apply the preceding lemma and theorem to prove that most of the complex reflection groups $G(r,p,n)$ with $\gcd(p,n) =2$ do not have generalized involution models.  We proceed by an argument similar to one used by Baddeley to prove that the Weyl group of type $D_{2n}$ does not have an involution model if $n>1$  \cite[Proposition 4.8.1]{B91}.  

\begin{lemma}\label{last2}
Let $r,p,n$ be positive integers with $p$ dividing $r$.  If $\gcd(p,n) = 2$, then $G(r,p,n)$ has a generalized involution model  if and only if $n=2$ and $r/p$ is odd.
\end{lemma}


\begin{proof}
Assume $\gcd(p,n) = 2$ so that $n$ and $r$ are both even.  Given Lemma \ref{last1},  we may assume that  $r/p$ is odd.  Suppose $\cM = \{ \lambda_i : H_i \to \CC\}$ is a generalized involution model for $G = G(r,p,n)$ with respect to some automorphism $\tau \in \Aut(G)$.  Then each $H_i = C_{G,\tau}(\omega_i)$ for a set of orbit representatives $\omega_i \in \cI_{G,\tau}$, and by Theorem \ref{thm4.5} and Lemmas \ref{classifying} and \ref{cmplx}, we may assume that $\tau$ is the usual inverse transpose automorphism $g\mapsto \overline g$.  Identify $\ZZ_2 \subset \ZZ_r$ as the subgroup $\ZZ_2 = \{0,r/2\} = (r/2)\ZZ_r$, so that we can view $G(2,2,n)$ as a subgroup of $G$.    

If we define $\Phi:  G(r,1,n) \to S_n$ as the surjective homomorphism given by $\Phi :(x,\pi) \mapsto \pi$, then $\Phi$ restricts to a surjective homomorphism $G(r,p,n) \to S_n$, and it follows from Theorem 5.2 in \cite{M}$-$also, it is not difficult to see directly$-$that the image under $\Phi$ of each $\tau$-twisted centralizer $H_i$ is conjugate to some subgroup $ W_k \times S_i \times S_j$ in $\sG(n)$.  
Thus the model $\cM_\Phi$ for $S_n$ defined by Theorem \ref{4.2.3} is based on $\sG(n)$.   

We now observe that  the generalized involutions 
\[ e = \( (0,\dots,0), 1\) = 1 \in\cI_{G,\tau}\qquad\text{and}\qquad \omega=\( (1,-1,1,-1\dots,1,-1), 1\) \in\cI_{G,\tau}\] belong to disjoint twisted conjugacy classes, since every element of the orbit of $e$ is of the form $(x,1) \in G$ with $x_i \in 2\ZZ_r$ for all $i$.   Since the stabilizers of elements in a given orbit are all conjugate, we may therefore assume without loss of generality that one linear character of $C_{G,\tau}(e)$ appears in $\cM$ and one linear character of $C_{G,\tau}(\omega)$ appears in $\cM$.  

  Because $r/p$ is odd, so that $r/2 \in \ZZ_2\subset \ZZ_r$ is an odd multiple of $p/2$, we have
\[ C_{G,\tau}(e)=(\ZZ_2 \wr S_n) \cap G = G(2,2,n).\]  To calculate $C_{G,\tau}(\omega)$, we observe that if 
$ z = \((1,0,1,0,\dots,1,0),1\) \in  G(r,1,n)$ and $c\in G(r,1,n)$ is the central element defined by (\ref{gen}),
then $z\cdot e \cdot\leftexp{\tau}z^{-1} = \omega c$.  Consequently, if $g \in G$ then $g\cdot\omega \cdot \leftexp{\tau}g^{-1} = \omega$ if and only if $g\cdot\omega c \cdot \leftexp{\tau}g^{-1} = \omega c$, and so $C_{G,\tau}(\omega) = \Ad(z) \(G(2,2,n)\)$. 

The group $C_{G,\tau}(e)=G(2,2,n)$ has only two linear characters $\lambda_1$ and $\lambda_2$, given by restricting the linear characters  $\One_{\ZZ_r} \wr (n)$ and $\One_{\ZZ_r} \wr (1^n)$ of $G(2,1,n) = \ZZ_2 \wr S_n$, respectively.  It is evident from the definition of these characters that $\ker (\Phi) \cap G(2,2,n)= (\ZZ_r)^n \cap G(2,2,n) \subset \ker(\lambda_i)$ for $i=1,2$ and that
\[ \cR_\Phi(\lambda_1) = \chi^{(n)} = \One \in \Irr(S_n)\qquad\text{and}\qquad \cR_\Phi(\lambda_2) = \chi^{(1^n)} = \sgn \in \Irr(S_n).\]   
Let $\lambda_i' = \lambda_i \circ \Ad(z)^{-1}$; then $\lambda_1',\lambda_2'$ are the only linear characters of $C_{G,\tau}(\omega)$, and since $\Phi \circ \Ad(z)^{-1} = \Phi$ as $z \in \ker(\Phi)$,  we have $\cR_\Phi(\lambda_i') = \cR_\Phi(\lambda_i)\neq 0$.  Thus either $\cR_\Phi(\lambda_1), \cR_\Phi(\lambda_2') \in \cM_\Phi$ or $\cR_\Phi(\lambda_2), \cR_\Phi(\lambda_1') \in \cM_\Phi$.  In particular $\One$ and $\sgn$ must both appear in $\cM_\Phi$, so by Lemma \ref{baddeley1} we have $n=2$.  In this case we know by Corollary \ref{thm-cor} that $G$ indeed has a generalized involution model, which completes the proof.
\end{proof}

\subsection{Classification of Groups with Generalized Involution Models}

We may now prove the theorem promised in the introduction.

\begin{theorem}\label{final}
A finite complex reflection group has a generalized involution model if and only if each of its irreducible factors is one of the following:
\begin{enumerate}
\item[(i)] $G(r,p,n)$ with $\gcd(p,n)=1$.
\item[(ii)] $G(r,p,2)$ with $r/p$ odd.
\item[(iii)] $G_{23}$, the Coxeter group of type $H_3$.
\end{enumerate}
\end{theorem}

\begin{proof}
Let $G$ be a finite complex reflection group.  Then $G$ is a product of irreducible complex reflection groups, so by Lemma \ref{products} it suffices to prove that the only irreducible complex reflection groups are those of types (i), (ii), and (iii).

To this end, first suppose $G =G_i$ for some $4\leq i \leq 37$ is an exceptional irreducible complex reflection group. If $G=G_{23}$ is the Coxeter group of type $H_3$, then $G$ has a generalized involution model by Corollary \ref{coxeter-cor}.  To prove no other exceptional groups have generalized involution models, we resort to an exhaustive computer search using the {\tt{GAP}} package {\tt{CHEVIE}}.  Several fortunate circumstances make this computation tractable.  First, by Corollary \ref{coxeter-cor}, we do not need to examine the Coxeter groups $G_{28}$, $G_{30}$, $G_{35}$, $G_{36}$, and $G_{37}$.  Second, it follows from Theorem \ref{bg-thm} that the exceptional groups $G_{27}$, $G_{29},$ and $G_{34}$ do not have generalized involution models because, upon examination of their character tables, one finds that if $G$ is one of these groups then $\sum_{\psi \in \Irr(G)} \psi$ assumes negative values.  Checking that each of the remaining exceptional groups does not have a generalized involution model by a brute force search is a feasible and not very time consuming calculation.  In particular, by Lemmas \ref{classifying} and \ref{cmplx} one only needs to examine at most one automorphism  for each group; we list candidates for this automorphism in Table \ref{tbl1}.  The remaining exceptional groups neither are prohibitively large nor have an excessive number of twisted conjugacy classes.

To deal with the infinite series, suppose $G = G(r,p,n)$ for some positive integers $r,p,n$ with $p$ dividing $r$.  If $\gcd(p,n)\leq 2$ then it follows from Corollary \ref{thm-cor} and Lemma \ref{last2} that $G$ has a generalized involution model if and only if $G$ is of the form (i) or (ii).  We may therefore assume $\gcd(p,n)>2$, so that $n>2$ and $r>2$.  

Suppose $G$ has a generalized involution model with respect to some  $\upsilon \in \Aut(G)$ with $\upsilon^2=1$.  By Theorem \ref{bg-thm} we then have $\sum_{\psi \in \Irr(G)} \psi(1) = | \cI_{G,\upsilon}|$ and $\epsilon_\upsilon(\psi)=1$ for all $\psi \in \Irr(G)$, so by Proposition 2 in \cite{BG2004} the elements $g^{-1}$ and $\leftexp{\upsilon}g$ are conjugate for all $g \in G$. It follows that $\upsilon$ preserves the normal subgroup $N=(\ZZ_r)^n \cap G$, and so by Lemma \ref{main-aut} we can write \[\upsilon = \Ad(g) \circ \alpha_{j,k,z},\qquad\text{ for some $g \in G(r,1,n)$ and $j,k,z$ as in (\ref{aut-prop-cond}).}\]  For some $a \in \ZZ_r$ we have $gt^{-a} \in G$, and if we let $\upsilon' = \Ad(t^a) \circ \alpha_{j,k,z}$ then  $g^{-1}$ and $\leftexp{\upsilon'}g$ are conjugate for all $g \in G$.  This fact implies that $\upsilon$ is the composition of an inner automorphism with the inverse transpose automorphism.

To see this, observe that $\Ad(t^a)$ fixes all element of $N$.  Therefore, if $x = \( e_1 -2e_2+e_3,1\) \in N$ then $\leftexp{\upsilon'} x = \( j(e_1-2e_2+e_3),1\)$, while all conjugates of $x^{-1}$ in $G$ are of the form $\( -e_{i_1} + 2e_{i_2} - e_{i_3},1\)$ for distinct $i_1,i_2,i_3 \in [1,n]$.  Since $r>2$, we must have $j\equiv -1 \modu r)$, and we may assume $j=-1$.  If $p=r$ then $\alpha_{j,k,z} = \alpha_{j,0,z}$ for all $k$.  If $p<r$, then 
\[ \leftexp{\upsilon'} t^p = c^{pk} t^{-p} = \( p(k-1) e_1 + pk (e_2 + \dots + e_n), 1\)\] while all conjugates of $t^{-p}$ are of the form $\( -p e_i, 1\)$ for $i \in [1,n]$.  Since $n>2$, it follows that $t^{-p}$ and $\leftexp{\upsilon'}t^p$ are conjugate only if $pk=0$ in $\ZZ_r$, in which case $\alpha_{j,k,z}  = \alpha_{-1,0,z}$.  As $\Ad(t^a)(s_2) = s_2$, we have $\leftexp{\upsilon'} s_2 = z s_2$ while all conjugates of $s_2^{-1}=s_2$ in $G$ are of the form $\( b e_{i_1} - b e_{i_2}, (i_1\cyc i_2) \)$ for $b \in \ZZ_r$ and $1\leq i_1 < i_2 \leq n$.  Again since $n>2$, it follows that $\leftexp{\upsilon'} s_2$ and $s_2^{-1}=s_2$ are conjugate only if $z = 1$.  Thus $\alpha_{j,k,z} = \alpha_{-1,0,1}$ is precisely the inverse transpose map.
Furthermore, since $\alpha_{-1,0,1}$ fixes all elements of $S_n$ and since each $\pi \in S_n$ is conjugate in $S_n$ to $\pi^{-1}$, it follows from Lemma \ref{inner} that $\Ad(t^a)$ defines an inner automorphism of $G$.  

We therefore  may assume that $\upsilon = \Ad(g) \circ \tau$ where $g \in G$ and $\tau :g\mapsto\overline g$ is the inverse transpose automorphism.  We now observe that $\omega \in G$ has $\omega \cdot \leftexp{\upsilon}\omega = 1$ if and only if $(\omega g) \cdot \leftexp{\tau} (\omega g) = g\cdot \leftexp{\tau} g$, so $|\cI_{G,\upsilon}|=
  | \{ \omega \in G : \omega \cdot \leftexp{\tau} \omega = g\cdot \leftexp{\tau} g\} |.$  
Since $\tau=\tau^{-1}$ and $\upsilon^2=1$, the element $g\cdot \leftexp{\tau} g$ is central, and 
as $n>2$, this implies that $|g\cdot \leftexp{\tau}g|=1$.  
Define $\cX_\pi(h) \subset (\ZZ_r)^n$ for each fixed $\pi \in S_n$ and $h \in G$ as the set of $x \in (\ZZ_r)^n$ with \[ x_1 + \dots + x_n \in p \ZZ_r\qquad\text{and}\qquad (x,\pi)\cdot\leftexp{\tau}(x,\pi)  = \( \pi^{-1} (x) - x, \pi^2\) = h\cdot\leftexp{\tau} h.\]
If $x,y \in \cX_\pi(g)$ then $x-y \in \cX_\pi(1)$ since  if $\cX_\pi(g)$ is nonempty then $\pi^2=|g\cdot \leftexp{\tau}g|=1$.  Hence $|\cX_\pi(g)| \leq |\cX_\pi(1)|$ for all $\pi \in S_n$.  Since
$\{ \omega \in G : \omega \cdot \leftexp{\tau}\omega = h\cdot\leftexp{\tau} h\} =\{ (x,\pi) : \pi \in S_n,\ x \in \cX_\pi(h)\},$ it follows that 
\[
|\cI_{G,\upsilon}|= | \{ \omega \in G : \omega \cdot \leftexp{\tau} \omega = g\cdot \leftexp{\tau} g\} | = \sum_{ \pi \in S_n} |\cX_\pi(g)| \leq \sum_{\pi \in S_n} |\cX_\pi(1)| = | \{ \omega \in G : \omega \cdot \leftexp{\tau} \omega = 1\}|.\]
We thus have $\sum_{\psi \in \Irr(G)} \psi(1) = |\cI_{G,\upsilon}| \leq |\cI_{G,\tau}|$. 
By Theorem \ref{thm4.5} this inequality must become equality, which contradicts the assumption that $\gcd(p,n)>2$.  We conclude that the only irreducible groups with generalized involution models are those of types (i)-(iii), which completes our proof. \end{proof}

We conclude with a few questions.  First, in \cite{C} Caselli 
defines $G(r,p,q,n)$ for positive integers $r,p,q,n$ with $p,q$ dividing $r$ and $pq$ dividing $rn$ to be the quotient group $G(r,p,n) / \langle c^{r/q} \rangle$.  He calls such groups \emph{projective complex reflection groups} and investigates their Gelfand models in \cite{C2009}.  The normal subgroup $\langle c^{r/q} \rangle \vartriangleleft G(r,p,n)$ is closed under transposes, so we have a well-defined notion of a transpose for elements of $G(r,p,q,n)$.  Among other results, Caselli proves that 
\[ \sum_{\psi \in \Irr\(G(r,p,q,n)\)} \psi(1) = | \{ \omega \in G(r,p,q,n) : \omega^T = \omega\}|\] if and only if $\gcd(p,n) \leq 2$, or $\gcd(p,n)=4$ and $r\equiv p \equiv q \equiv n \equiv 4 \modu 8)$.  This begs the following question.

\begin{question} In which of these cases does $G(r,p,q,n)$ have a generalized involution model with respect to the inverse transpose automorphism?  
 More broadly, which projective complex reflection groups have generalized involution models?  
\end{question}

Second, in \cite{APR2007} Adin, Postnikov, and Roichman describe 
a representation $\rho_q$ for the one-parameter Iwahori-Hecke algebra of the symmetric group which is a Gelfand model when the algebra is semisimple and which becomes the Gelfand model
$\rho_{1,1,n}$ for $S_n = G(1,1,n)$ when one specializes to $q=1$.  
We wonder if this situation can be extended to other complex reflection groups, and echoing a question posed at the end of \cite{APR2008}, ask the following.


\begin{question} Can one describe a similar deformation of the representations $\rho_{r,p,n}$ and $\wt \rho_{r,p,n}$ or of the representations appearing in \cite{C2009} to obtain a Gelfand model for the Iwahori-Hecke algebras of the irreducible Weyl groups, or for the  Ariki-Koike algebras of $G(r,p,n)$ when $\gcd(p,n)\leq 2$?
\end{question}

\end{document}